\documentclass[12pt]{article}

\usepackage{latexsym}
\usepackage{calc}
\usepackage{amssymb, amsmath, amsfonts, listings, mathrsfs} 
\usepackage{graphicx}
\usepackage{float}
\usepackage[pagebackref]{hyperref}
\usepackage{booktabs}
\usepackage{threeparttable}
\usepackage{makecell}
\usepackage{multirow}
\usepackage{algorithm}
\usepackage[noend]{algpseudocode}

\hypersetup{
  colorlinks=true,
  linkcolor=blue,
  filecolor=magenta,
  urlcolor=cyan,
  citecolor=blue
}

\renewcommand*{\backrefalt}[4]{%
    \ifcase #1 \footnotesize{(Not cited.)}%
    \or        \footnotesize{(Cited on page~#2)}%
    \else      \footnotesize{(Cited on pages~#2)}%
    \fi}

\usepackage{lmodern}
\usepackage{geometry}
\newcounter{hours}\newcounter{minutes}

\def\argmin{\mathop{\rm argmin}}

\def\argmax{\mathop{\rm argmax}}

\def\nr{\par \noindent}

\def\Def{\stackrel{\mathrm{def}}{=}}

\def\sign{{\rm sign \,}}

\def\dom{{\rm dom \,}}
\def\beq{\begin{equation}}
\def\eeq{\end{equation}}

\def\R{\mathbb{R}}
\def\E{\mathbb{E}}

\def\BI{\begin{itemize}}
\def\EI{\end{itemize}}

\newtheorem{theorem}{Theorem}
\newtheorem{lemma}{Lemma}
\newtheorem{corollary}{Corollary}

\newtheorem{proposition}{Proposition}
\newtheorem{assumption}{Assumption}
\newtheorem{definition}{Definition}

\newtheorem{example}{Example}
\newtheorem{remark}{Remark}
\newcommand{\proof}{\bf Proof: \rm \nr}
\newcommand{\qed}{\hfill $\Box$ \nr \medskip}

\renewcommand\arraystretch{1}
\def\ba{\begin{array}}
\def\ea{\end{array}}
\def\beann{\begin{eqnarray*}}
\def\eeann{\end{eqnarray*}}
\def\bea{\begin{eqnarray}}
\def\eea{\end{eqnarray}}

\def\BT{\begin{theorem}}
\def\ET{\end{theorem}}
\def\BL{\begin{lemma}}
\def\EL{\end{lemma}}
\def\BC{\begin{corollary}}
\def\EC{\end{corollary}}
\def\BE{\begin{example}}
\def\EE{\end{example}}
\def\BD{\begin{definition}}
\def\ED{\end{definition}}
\def\BR{\begin{remark}}
\def\ER{\end{remark}}
\def\BAS{\begin{assumption}}
\def\EAS{\end{assumption}}
\def\BI{\begin{itemize}}
\def\EI{\end{itemize}}
\def\BP{\begin{proposition}}
\def\EP{\end{proposition}}

\def\BMP{\begin{minipage}{9.5cm}}
\def\EMP{\end{minipage}}
\def\MPT{\begin{minipage}{11.5cm}}
\def\EPT{\end{minipage}}

\def\la{\langle}
\def\ra{\rangle}

\def\QF{\hspace{5ex} \Box}
\def\QR{\hfill \Box}

\usepackage{color}
\usepackage{pifont}
\definecolor{mydarkgreen}{RGB}{39,130,67}

\definecolor{mydarkred}{RGB}{192,47,25}

\algnewcommand\algorithmicon{\textbf{for}}
\algnewcommand\algorithmicfrom{\textbf{from}}
\algblockdefx[ON]{On}{EndOn}[1]
  {\algorithmicon\ #1\ \ \algorithmicdo}
  {\algorithmicuntil\ }

\newcommand{\cmark}{{\normalsize {\color{mydarkgreen}\ding{51}}}}%
\newcommand{\xmark}{{\normalsize {\color{mydarkred} \ding{55}}}}%

\title{
\textbf{Super-Universal Regularized\\ Newton Method\thanks{This project has received funding from the European Research Council
		(ERC) under the European Union's Horizon 2020 research and innovation
		programme (grant agreement No. 788368).
	The second author acknowledges support by the French government under the management of the Agence Nationale de la Recherche as part of the ``Investissements d’avenir'' program, reference ANR-19-P3IA-0001 (PRAIRIE 3IA Institute).
	The research of the third author was also supported by Multidisciplinary Institute in Artificial intelligence MIAI@Grenoble Alpes (ANR-19-P3IA-0003).}}
}

\author{Nikita Doikov\thanks{ICTEAM, UCLouvain, Louvain-la-Neuve, Belgium.
	 	nikita.doikov@uclouvain.be} 
		\and 
		Konstantin Mishchenko\thanks{CNRS, ENS, Inria Sierra, Paris, France. konsta.mish@gmail.com} 
		\and Yurii Nesterov\thanks{CORE, UCLouvain, Louvain-la-Neuve, Belgium.
			yurii.neterov@uclouvain.be}
		}

\date{ August 11, 2022 }

\begin{document}
\maketitle

\begin{abstract}
	We analyze the performance of a variant of
	Newton method with quad\-ratic regularization
	for solving composite convex minimization problems.
	At each step of our method, we choose regularization parameter
	proportional to a certain power of the gradient norm at the current point.
	We introduce a family of problem classes characterized by
	H\"older continuity of either the second or third derivative.
	Then we present the method with a simple adaptive search procedure allowing an automatic adjustment to the problem class with the best global complexity bounds, without knowing specific parameters of the problem.
	In particular, for the class of functions with Lipschitz continuous third derivative,
	we get the global $O(1/k^3)$ rate, which was previously attributed to
	third-order tensor methods. 
	When the objective function is uniformly convex, we justify
	an automatic acceleration of our scheme, resulting
	in a faster global rate and local superlinear convergence.
	The switching between the different rates (sublinear, linear, and superlinear)
	is automatic. Again, for that, no a priori knowledge of parameters
	is needed.
\end{abstract}

\textbf{Keywords:} Newton method, regularization, global convergence, convex optimization,
global complexity bounds, universal methods

\newpage
\section{Introduction}

\paragraph{Motivation.}

Newton's method is one of the most important tools 
in Numerical Analysis and Continuous Optimization.
It has a reputation for being a powerful algorithm,
especially due to its ability 
to solve ill-conditioned problems.
The method has a local quadratic convergence,
thus converging extremely fast
in a neighbourhood of the solution \cite{kantorovich1948functional}.
However, the global behaviour of Newton's method 
has been remaining an active area of research
for several decades.

It is widely known that
the classical Newton method with a unit stepsize
may not converge globally, even if the problem 
is strongly convex
(see, e.g., Example 1.4.3 in \cite{doikov2021new}).
Consequently, there were many techniques developed
for the method to improve its global behaviour, 
including \textit{damped} Newton steps
combined with line search strategies~\cite{kantorovich1948functional,ortega2000iterative},
Levenberg-Marquardt regularization~\cite{levenberg1944method,marquardt1963algorithm},
and \textit{trust-region} approach~\cite{goldfeld1966maximization,conn2000trust}.
(See also \cite{polyak2007newton} for an extensive historical overview.)
However, it was still difficult to establish global complexity guarantees
that are provably better than that of the Gradient Methods.

A major shift in the paradigm has been made after the work \cite{nesterov2006cubic},
where \textit{cubic regularization} of Newton's method (CNM) with its global convergence guarantees
was developed. The main idea was to start with a particular problem class,
the functions with Lipschitz continuous Hessian,
which naturally leads to a globally convergent second-order scheme.
The subproblem becomes the minimization of a quadratic model of the function
augmented by the third power of Euclidean norm.
While each iteration of the method requires solving a univariate nonlinear equation, 
the arithmetical cost of such an operation remains of the same order
as for the standard Newton step.

Later on, \textit{adaptive}~\cite{cartis2011adaptive1,cartis2011adaptive2} 
and \textit{universal}~\cite{grapiglia2017regularized,doikov2021minimizing}
second-order methods based on cubic regularization with 
an adjustment of the Lipschitz constant were developed.
In \cite{grapiglia2017regularized}, it was shown that
the adaptive search makes the CNM work properly on functions 
with H\"older continuous Hessian, automatically achieving 
the correct global complexity, and in \cite{doikov2021minimizing}
the universality of CNM was studied on uniformly convex functions.

A parallel line of work was done for the Newton method with quadratic regularization.
In \cite{polyak2009regularized}, the author proposed to use
the gradient norm as a regularization coefficient, which preserves
the local quadratic convergence of the Newton iterations.
However, to ensure a global rate, it was needed to use some
damped steps, which make the convergence slower
than that of CNM.
The idea to approximate the cubic step by a quadratic regularization
probably appeared for the first time in \cite{ueda2009regularized},
still having a worse rate. Eventually, it was first proven in \cite{mishchenko2021regularized},
and independently rediscovered in \cite{doikov2021gradient},
that the use of \textit{square root of the gradient norm}
as the regularization coefficient provides the method with 
the \textit{fast global rate} of CNM, while each iteration
requires now just one standard matrix inversion.

Another emerging trend in Optimization 
has been to use higher-order Taylor's model of the objective,
which potentially would result in even more powerful methods,
called \textit{Tensor Methods} (\cite{baes2009estimate,birgin2017worst,nesterov2019implementable,grapiglia2020tensor,dvurechensky2019near,cartis2020sharp}).
The price for such an advancement is clear:
the subproblem, which is a minimization of the 
high-order polynomials, becomes more and more difficult.
A valuable observation was made in \cite{nesterov2019implementable},
showing that the regularized Taylor polynomial 
of a convex function is convex,
which makes the subproblem solvable.
An efficient procedure for computing the third-order tensor step
was also proposed there. Following this direction, 
and utilizing the fact that \textit{the third derivative of a convex function is weak},
there were developed efficient third-order type schemes
which use only the \textit{second-order information}
\cite{nesterov2021superfast,kamzolov2020near,nesterov2022quartic}.

In this paper, we develop a surprisingly simple but very powerful
regularization strategy for Newton's method,
that provides the method with provably fast and universal global convergence rates.

The main idea behind our algorithm is to regularize
the second-order model of the objective by
the square of Euclidean norm, with regularization 
coefficient proportional to a certain power of the gradient norm.
In the simplest case of unconstrained minimization $\min_{x \in \R^n} f(x)$,
one iteration of our method is as follows:
$$
\boxed{
	\ba{rcl}
	\lambda_k & = & H_k \|\nabla f(x_k)\|^{\alpha}, \quad
	x_{k + 1} \; = \; x_k 
	- \bigl( \nabla^2 f(x_k) +  \lambda_k I   \bigr)^{-1} \nabla f(x_k), \quad k \geq 0,
	\ea
}
$$
where $\nabla f(x_k)$ is the current gradient, $\nabla^2 f(x_k)$ is the current Hessian, and $I$ is the identity matrix. In this scheme,
the power $\alpha$ can be fixed arbitrarily from the range $[\frac{2}{3}, 1]$. At the same time, the regularization constant $H_k$ is adjusted automatically by a standard backtracking procedure, based on the following stopping criterion:
$$
\boxed{
\ba{rcl}
\la \nabla f(x_{k+1}), x_{k} - x_{k + 1} \ra & \geq & {1 \over 4 \lambda_k} \| \nabla f(x_{k+1}) \|^2.
\ea
}
$$
Thus, at each iteration, the method needs to compute the Hessian once,
and the average number of the search steps is only \textit{two}.
The algorithm uses matrix inversion as the basic subroutine,
which can be implemented either with Linear Algebra tools, or by using 
a gradient-type solver in the large-scale setting.

We show that our strategy works for a wide range of problem classes,
characterized by H\"older continuity of \textit{either}
the third \textit{or} second derivative.
The algorithm itself does not need to know any parameters of the problem class.
Therefore, our method
\textit{automatically} achieves the rates of convergence of the
Gradient Method, Cubic Newton,
and third-order Tensor Methods
on the corresponding problem classes.
Moreover, when the objective is strictly convex,
our new algorithm makes sure to get an acceleration,
automatically switching between superlinear, linear, and superlinear rates.
We attribute the name \textit{Super-Universal} to 
a method possessing all these features.

\subsection{Contributions} 

\begin{table}[]
\centering
\setlength\tabcolsep{4pt}
  \begin{threeparttable}[b]
    \centering
    \caption{A conceptual comparison of our method to basic deterministic algorithms. 
    Since accelerated methods work on the same problem classes, we do not include them here. 
    For the convergence rates of our method
    on convex functions, see Corollary~\ref{SummaryCor}.}
    \label{tab:comparison}
    {
    \footnotesize
    \begin{tabular}{cccccccc}
    \toprule
		\multirow[c]{2}{*}{Method} & 
		\hspace*{-5pt}\multirow[c]{2}{*}{\makecell{Easy\\ implementation}} & \multirow[c]{2}{*}{\makecell{Local superlinear\\ convergence}} & \multicolumn{3}{c}{Global convergence} & 
		\hspace*{-5pt}\multirow[c]{2}{*}{References} \\ \cmidrule(lr){4-6}
		&  &  & \multicolumn{1}{c|}{Lip.\ grad.} & \multicolumn{1}{c|}{Lip.\ Hess.} & Lip.\ 3rd der. \\ \midrule
		Gradient Method & \cmark & \xmark & \cmark & \xmark & \xmark & \cite{nesterov2018lectures} \\
		Classical Newton & \cmark & \cmark & \xmark & \xmark & \xmark & \cite{nesterov2018lectures} \\ 
		Universal Cubic Newton & \xmark & \cmark & \cmark & \cmark & \xmark & \cite{nesterov2006cubic,grapiglia2017regularized}  \\ 
		\makecell{
			Universal 3rd-order \\
			Tensor Method} & \xmark & \cmark & \xmark & \cmark & \cmark & \cite{grapiglia2020tensor,doikov2021local}  \\ \midrule
		Super-Universal Newton & \cmark & \cmark & \cmark & \cmark & \cmark & \textbf{Ours} \\
        \bottomrule
    \end{tabular}
    }
  \end{threeparttable}
\end{table}  

We present our theory gradually starting with basic results and eventually leading to Super-Universal Method. We start with Section~\ref{SectionRegNewton}, where we introduce the composite optimization problem
and discuss properties of the regularized Newton method.

In Section~\ref{SectionProblemClasses}, we present a family of problem classes 
characterized by H\"older continuity of the second and third derivatives of the 
smooth part of the objective.
We provide a univariate parametrization $2 \leq q \leq 4$ of these classes, introducing 
the corresponding smoothness parameter $M_q$.
The particular cases include
Lipschitz continuity of the third derivative $(q = 4)$,
Lipschitz continuity of the Hessian $(q = 3)$,
and bounded variation (or boundedness) of the Hessian $(q = 2)$.

Section~\ref{SectionGradReg} contains our main tool,
that is the choice of regularization coefficient proportional to a certain power
of the gradient norm. We prove several inequalities leading to the global
and local convergence of the basic steps,
which lead to a simple iterative scheme given in Algorithm~\ref{alg:FixedNewton}.

In Section~\ref{SectionUniversal},
we develop Super-Universal Newton Method (Algorithm~\ref{alg:AdaptiveNewton})
with an adaptive search procedure based on a new stopping
criterion. 
Our method does not need to know any parameters of the objective,
and achieves a \textit{universal} global complexity for all our problem classes.
For general convex case,
the convergence rate is $O(k^{1 - q})$
in terms of the functional residual,
where $k$ is the iteration counter.

In Section~\ref{SectionSConvex},
we study the global and local convergence of our method on subclasses
of strictly convex functions.
We introduce a new characteristic of optimization problems 
called $s$-\textit{relative size}
$D_s$ ($s \geq 2$).
Our definitions clarify the standard notion of uniform convexity
and allow continuous change in the convexity degree.
For $s = 2$, our assumption implies strong convexity, and
for $s = \infty$ it means that the initial level set
is bounded with diameter $D \equiv D_{\infty}$.
We show that our method achieves automatically improved rates on these subclasses.

The following table provides a summary of the global complexity guarantees.
We are interested in the number of iterations to reach $\varepsilon$
accuracy in terms of the functional residual
(presenting only the main terms and omitting absolute constants).

\begin{table}[!h]
	\caption{\small Global Complexity of  Super-Universal Newton Method (Algorithm~\ref{alg:AdaptiveNewton}). If $s=2$, it means the objective is strongly convex, and we have superlinear convergence. $V_F$ denotes 
the size of the initial level set measured by symmetrized 
Bregman divergence, defined by \eqref{VFDef}.}
	\begin{center}
		{
			\small
			\centering
			\renewcommand{\arraystretch}{1.5}
			\begin{tabular}{ | c | c | c | c | }
				\hline
				$\;\; 2 \leq s < q \;\;$ & 
				$\;\;\;\; s = q \;\;\;\;$ & 
				$\;\; q < s < \infty \;\;$ &
				$ s = \infty $
				\\
				\hline
				\begin{tabular}{@{}l@{}}
					\\[-13pt]			
					$ \! \Bigl(\! M_q  \frac{D_s^s D^{q - s} }{V_F}\Bigr)^{\frac{1}{q - 1}} + \ln\ln\frac{1}{\varepsilon}$\!  \\[10pt]
				\end{tabular}
				&
				\begin{tabular}{@{}l@{}}
					\\[-13pt]			
					$ \Bigl(\! M_q \frac{D_q^q}{V_F} \Bigr)^{\frac{1}{q - 1}} \ln\frac{1}{\varepsilon}$  \\[10pt]
				\end{tabular}
				&
				$
				\begin{tabular}{@{}l@{}}
				\\[-13pt]			
				$ \Bigl(\! M_q \frac{ D_s^q }{ (V_F^{q} \varepsilon^{s - q})^{1/s} } \Bigr)^{\frac{1}{q - 1}}$\!\! \\[10pt]
				\end{tabular}
				$
				&
				$
				\begin{tabular}{@{}l@{}}
				\\[-13pt]			
				$ \! \Bigl( \! M_q  \frac{D^q}{\varepsilon} \Bigr)^{\frac{1}{q - 1}}$\!\!  \\[10pt]
				\end{tabular}
				$ 
				\\
				\hline
				
			\end{tabular}
		}
	\end{center}
\end{table}

\begin{figure}[H]
	\centering
	\includegraphics[width=0.6\textwidth ]{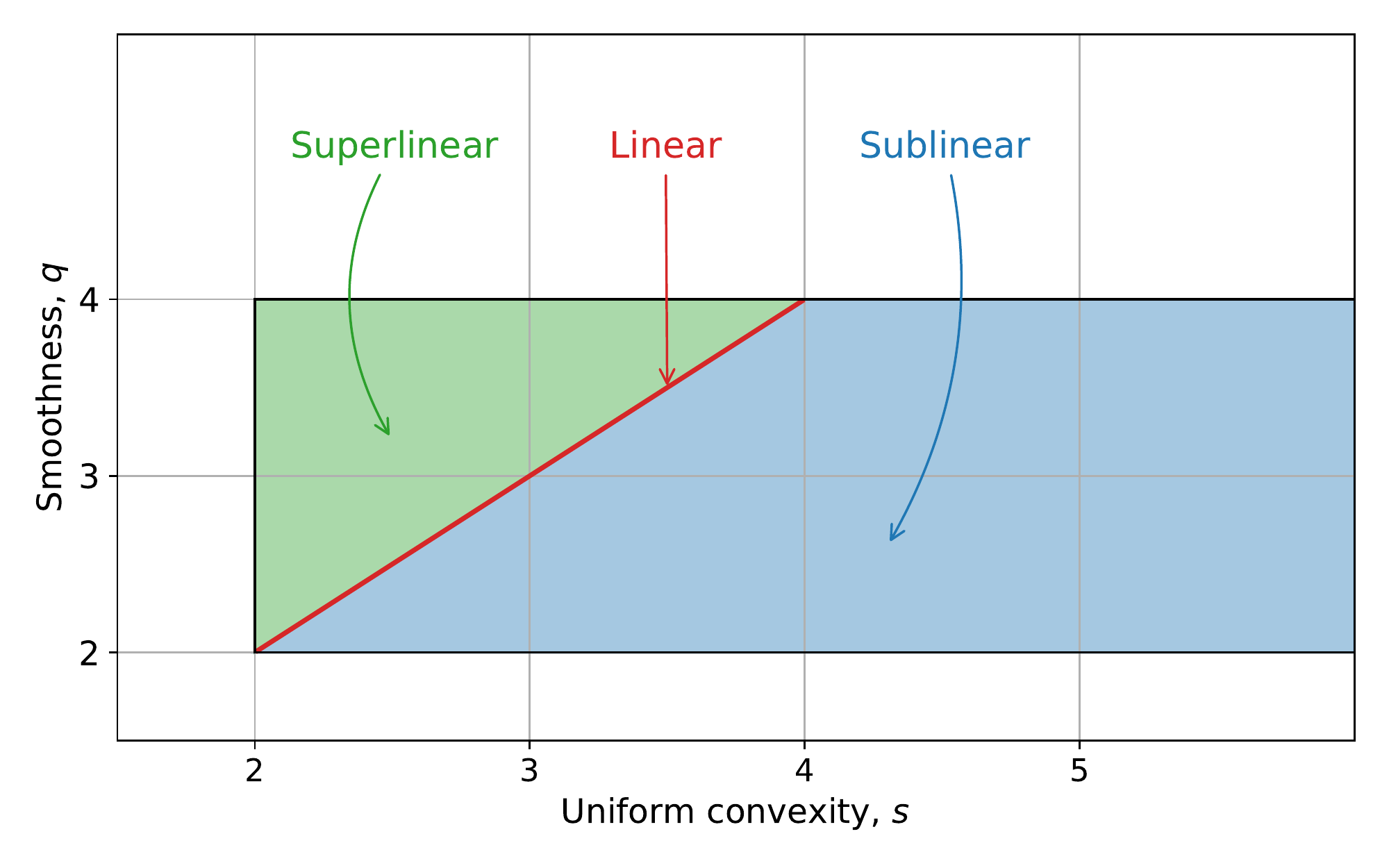}\\[2pt]
	\caption{ \small Global Convergence Diagram for Super-Universal Newton Method (Algorithm~\ref{alg:AdaptiveNewton}).}
	\label{FigureRegions}
\end{figure}

Qualitatively, the rates are split into two regions with a switching line in the rate of convergence given by the case $s = q$.
For $2 \leq s \leq q$, the complexity depends on target accuracy $\varepsilon$
only logarithmically.
When $s < q$, the method has a superlinear convergence,
and the case $s = q$ gives us the global linear rate.
For $s > q$, the dependence on accuracy is polynomial,
which corresponds to sublinear rates. 
The whole picture becomes two-dimensional taking into account the range for
the degree of smoothness: $2 \leq q \leq 4$ (see Figure~\ref{FigureRegions}).

\bigskip

Note that our new algorithm provides us with the worst-case complexity bounds for all known problem classes with computable second derivative supported till now by different 
first-, second- and third-order schemes.
Moreover, our Super-Universal Newton Method seems to eliminate the need of using 
non-accelerated third-order Tensor Methods
in Convex Optimization.

Finally, we present our numerical experiments in 
Section~\ref{SectionExperiments},
and provide a discussion on possible future developments
in Section~\ref{SectionDiscussion}.

\subsection{Notation}
In what follows, we denote by $\E$ a finite dimensional real vector space,
and by $\E^{*}$ its dual, which is the space of linear functions on $\E$.
The value of function $s \in \E^{*}$ on vector $x \in \E$ is denoted by
$\la s, x \ra$.

Let us fix a self-adjoint positive-definite linear operator $B\colon \E \to \E^{*}$ 
(notation $B = B^{*} \succ 0$), and use it to define the Euclidean norm for the primal space:
$$
\ba{rcl}
\| x \| & \Def & \la Bx, x \ra^{1/2}, \qquad x \in \E.
\ea
$$
Then, in the dual space we apply the induced norm:
$$
\ba{rcl}
\| s \|_{*} & \Def & \max\limits_{h \in \E \, : \, \|h\| \leq 1 }
\la s, h \ra
\;\; = \;\; \la s, B^{-1} s \ra^{1/2},
\qquad s \in \E^{*}.
\ea
$$

For convex function $f(\cdot)$, we use notation $\partial f(x)$ 
for its subdifferential at point $x \in \dom f$.
If $f \colon \dom f \to \R$
is several times differentiable,
we denote its gradient by $\nabla f(x)$ and its Hessian by $\nabla^2 f(x)$.
Note that
$$
\ba{rcl}
\nabla f(x) & \in & \E^{*}, \qquad 
\nabla^2 f(x)h \;\; \in \;\; \E^{*}, \qquad x \in \dom f \subseteq \E, \; h \in \E.
\ea
$$
Along a fixed direction $h \in \E$, we use the following notation
for the second directional derivative:
$$
\ba{rcl}
\nabla^2 f(x)[h]^2 & \Def & \la \nabla^2 f(x)h, h \ra \;\; \in \;\; \R, \qquad h \in \E. 
\ea
$$
The third derivative, which is a trilinear symmetric form, is denoted by $\nabla^3 f(x)$.
Thus, 
$$
\ba{rcl}
\nabla^3 f(x)[h_1, h_2, h_3] & \in & \R, \qquad h_1, h_2, h_3 \in \E.
\ea
$$
When $h = h_1 = h_2 = h_3$, 
we use a shorthand: $\nabla^3 f(x)[h]^3 \Def \nabla^3 f(x)[h, h, h]$.

For symmetric multilinear forms, we define the induced norms 
in the standard way. For example,
$$
\ba{rcl}
\| \nabla^2 f(x) \| & \Def & 
\max\limits_{\substack{h_1, h_2 \in \E \\ \|h_1\| \leq 1, \|h_2\| \leq 1}}
\la \nabla^2 f(x) h_1, h_2 \ra
\;\; = \;\; 
\max\limits_{h \in \E  \, : \, \|h\| \leq 1} \| \nabla^2 f(x)h \|_{*},
\ea
$$
and
$$
\ba{rcl}
\| \nabla^3 f(x) \| & \Def & 
\max\limits_{\substack{ h_1, h_2, h_3 \in \E \\
		\|h_i\| \leq 1, \forall i}}
	\nabla^3 f(x)[h_1, h_2, h_3]
\;\; = \;\;
\max\limits_{h \in \E \, : \, \|h\| \leq 1}
| \nabla^3 f(x)[h]^3|,
\ea
$$
where the last equation is proved in Appendix 1 of \cite{nesterov1994interior}.

\subsection{Overview of the Main Ideas}
Before we proceed to formal proofs and detailed explanation of the main results, let us first sketch them here to provide a high-level intuition. For simplicity, we start with a discussion on how to solve 
the problem $\min_x f(x)$ without any nonsmooth components.

\paragraph{Removing Third Derivatives.}
An early observation was made in \cite{nesterov2019implementable} 
that a third-order Tensor Method can be implemented 
using second-order oracle calls with an auxiliary 
procedure that computes the action $\nabla^3 f(x)[h]^2$ of the tensor of the third derivative to an arbitary vector $h$. 
This is in a way similar to how second-order Newton method can be implemented by running a first-order method on a quadratic subproblem. Where the similarity disappears, however, 
is that third-order subproblem requires only \textit{near-constant} number of inner iterations, 
independently of any other function properties.

This observation was also used in \cite{nesterov2021superfast} to design inexact third-order methods that rely 
\textit{solely on second-order oracle}, approximating the action of the third derivative by a finite difference.
 Following upon these results, we find it natural to ask: Is it possible to skip formulating third-order subproblems and show instead that a simple second-order method is sufficient to achieve faster convergence?

It turns out that the answer is positive. The first step is to notice that, as stated in Lemma~3 in~\cite{nesterov2019implementable}, when the third derivative of a convex function is Lipschitz continuous (with constant $M_4 > 0$), one can show that its impact is always bounded as follows:
\[
\ba{rcl}
	\nabla^3 f(x)[h]^3
	& \le &
	\frac{1}{\tau}\nabla^2 f(x)[h]^2 + \frac{\tau}{2} M_4 \|h\|^4 \qquad \text{for any } \tau > 0 \text{ and } x,h\in\E.
\ea
\]
In other words, the third-order term is fully controlled by a combination of second-order and fourth-order terms. We can immediately plug this in the global upper bound on function,
\[
\ba{rcl}
	f(y)
    & \le & 
    \underbrace{\textstyle f(x) + \langle \nabla f(x), y-x\rangle + \frac{1}{2}\nabla^2 f(x)[y-x]^2 + \frac{1}{6}\nabla^3 f(x)[y-x]^3}_{\textrm{Taylor approximation of }f(y)} + \frac{1}{24} M_4\|y-x\|^4,
\ea
\]
which holds for any $x, y,\in \E$, and get
\begin{equation}
\ba{rcl}
	f(y)
    & \le & 
    f(x) + \langle\nabla f(x), y-x\rangle + \bigl( \frac{1}{2} + \frac{1}{6\tau} \bigr)\nabla^2 f(x)[y-x]^2 
    + \frac{1 + 2\tau}{24} M_4 \|y-x\|^4. 
    \label{eq:second_order_upper_bound}
\ea
\end{equation}
Optimizing this upper bound exactly would give a Newton-type iteration, but it still has two disadvantages. Firstly, the coefficient in front of the Hessian is not $\frac{1}{2}$. At the same time,
the factor $\frac{1}{2}$
is the best from the local perspective since it is responsible for the local \textit{superlinear} convergence of the method.
Secondly, the last term in the upper bound makes it polynomial of degree four and gives us 
a subproblem which is not trivially to solve. 
In our work, we provide a solution to both of these challenges by showing that we do not actually need to minimize the upper bound exactly. Instead, we prove that a regularized Newton iteration is sufficient to decrease the functional values despite not being an exact minimizer of this upper bound.

\paragraph{Gradient Regularization.}
Our idea to simplify the upper bound in~\eqref{eq:second_order_upper_bound} stems from 
the prior results on simplifying the Cubic Newton iteration. In CNM, each iteration is obtained by solving
\[
\ba{c}
	\argmin\limits_{x} \Bigl\{
	\la \nabla f(x_k), x - x_k \ra + \frac{1}{2} \nabla^2 f(x_k)[x - x_k]^2 + \frac{M_3}{6}\|x - x_k\|^3
 \Bigr\}    \qquad \textrm{(Cubic Newton Step)} 
\ea
\]
with some constant $M_3>0$. Because of the cubic term, this equation has no closed-form solution. At the same time, it was shown in \cite{mishchenko2021regularized} and later in \cite{doikov2021gradient} that if we replace the cubic regularization
by a quadratic with an appropriately chosen coefficient, the method would still converge with the same rate. In particular, if we set $\lambda_k=\sqrt{\frac{M_3}{3}\|\nabla f(x_k)\|_{*}}$ and produce the iterates by solving 
\[
\ba{c}
	\argmin\limits_{x} \Bigl\{
	\la \nabla f(x_k), x - x_k \ra + \frac{1}{2} \nabla^2 f(x_k)[x - x_k]^2 + \frac{\lambda_k}{2}\|x - x_k\|^2
 \Bigr\}    \qquad \textrm{(Regularized Newton Step)}
\ea
\]
then the rate of convergence remains the same as for CNM. 
Conceptually, there is little difference between cubic and quartic regularization as it appears in~\eqref{eq:second_order_upper_bound}. Therefore, we can apply the same ideas and replace the fourth power by the second again. 
A simple derivation shows that to lift the quartic regularization to the quadratic one, we 
need to use $\lambda_k = H_k\|\nabla f(x_k)\|^\frac{2}{3}_{*}$ with a sufficiently large constant $H_k>0$. 
An immediate drawback is that this approach requires knowing the problem class to choose the power of the gradient norm in the expression for $\lambda_k$. Indeed, we have to choose whether to use $\lambda_k\propto \|\nabla f(x_k)\|^\frac{1}{2}_{*}$ or $\lambda_k\propto \|\nabla f(x_k)\|^\frac{2}{3}_{*}$, and this choice is not trivial, since our objective
can belong to several problem classes simultaneously.
Thus, we must ask: Is it possible to design a method that would not require to know the parameters of the problem?

\paragraph{Super Universality.} 
As we have stated, regularized Newton is applicable to several problem classes at the same time, but it requires different strategies for choosing $\lambda_k$. 

In prior works, however, it was shown that 
one regularization strategy can sometimes work properly for different problem classes.
In particular, \cite{grapiglia2017regularized} and \cite{doikov2021minimizing} studied CNM for functions whose Hessian is H\"older continuous, and the method does not need to know the H\"older parameter. One of our goals, therefore, is to make our Newton method adaptive not just to the class of functions based on whether it is second or third derivative that is Lipschitz continuous, but also to the H\"older constant within each of the classes.

By working in this direction, we managed to obtain a universal regularization rule that works across all considered function classes.
It turned out that we can use $\lambda_k = H_k\|\nabla f(x_k)\|_{*}^\alpha$ with \textit{any} $\alpha\in[\frac{2}{3}, 1]$
fixed in advance, without even knowing what properties the minimized function has. Parameter $H_k$ is adjusted 
automatically by a standard adaptive procedure. This approach  
makes our Newton method much more universal than any other prior regularization techniques,
as for solving unconstrained convex minimization problems.

\paragraph{Composite Optimization.} Finally, it is important in many applications to support additional 
nonsmooth components, such as constraints or $\ell_1$-regularization. This corresponds to minimizing $f(x)+\psi(x)$, where $\psi(x)$ is a nonsmooth function. In first-order optimization, one can use so-called proximal operator that may even have a closed-form solution. In second-order optimization, the subproblem becomes more difficult, 
which requires minimization of quadratic function with the extra component $\psi$. 
Thus, for the complex $\psi$, we can no longer solve the iteration as a simple linear system.
However, we can use first-order gradient-based solvers for computing an inexact iteration. 
The composite formulation covers more applications, and we still do not have to worry about which problem class function $f$ belongs to. This extension is more straightforward than the other steps as it has been already considered for regularized Newton method in~\cite{doikov2021gradient}.

\section{Regularized Newton Step}
\label{SectionRegNewton}

Consider the following composite optimization problem:
\beq \label{MainProblem}
\ba{rcl}
\min\limits_{x \in \dom \psi} \Bigl \{   F(x) & = & f(x) + \psi(x) \Bigr\}.
\ea
\eeq
where function $f(\cdot)$ is convex and several times differentiable, 
and $\psi \colon \E \to \R \cup \{ +\infty \}$ is a proper closed convex function
with $\dom \psi \subseteq \E$.
For some $\lambda > 0$, 
consider the step of a variant of Newton method with quadratic regularization:
$$
\ba{rcl}
T_{\lambda}(x) & \Def & \argmin\limits_{y} \Bigl\{
\la \nabla f(x), y - x \ra + \frac{1}{2} \nabla^2 f(x)[y - x]^2 + \frac{\lambda}{2}\|y - x\|^2 + \psi(y)
 \Bigr\},
\ea
$$

If the composite part is absent ($\psi(y) \equiv 0$), this step can be rewritten in 
an explicit form:
$$
\ba{rcl}
T_{\lambda}(x) & = & x - \bigl( \nabla^2 f(x) + \lambda B  \bigr)^{-1} \nabla f(x),
\ea
$$
which is often called the Levenberg-Marquardt regularization 
\cite{levenberg1944method,marquardt1963algorithm}.
In the presence of $\psi(\cdot)$,
the point $T = T_{\lambda}(x)$ satisfies the following stationary condition:
\beq \label{StatCond}
\ba{rcl}
\la \nabla f(x) + \nabla^2 f(x)(T - x) + \lambda B(T - x), y - T \ra + \psi(y) & \geq & 
\psi(T),
\ea
\eeq
for any $y \in \dom \psi$. In other words, 
\beq \label{SubgDef}
\ba{rcl}
\psi'(T) & \Def & - \nabla f(x)  - \nabla^2 f(x)(T - x) - \lambda B(T - x) 
\;\; \in \;\; \partial \psi(T),
\ea
\eeq
and therefore,
\beq \label{def-DF}
\ba{rcl}
F'(T) & \Def & \nabla f(T) + \psi'(T)  \;\; \in \;\; \partial F(T).
\ea
\eeq
We also denote $r \Def \|T - x\|$.

Let us derive some inequalities for one step of the method.
Let $\mu \geq 0$ be a uniform bound for the minimal eigenvalue of the Hessian: 
$\nabla^2 f(x) \succeq \mu B$, $\forall x \in \dom \psi$. 
Note that, for any $s \in \partial \psi(x)$, it holds
\beq \label{HessBound1}
\ba{rcl}
\la \nabla^2 f(x)(T - x), T - x \ra
& \overset{\eqref{SubgDef}}{=} &
\la \nabla f(x) + \psi'(T), x - T \ra - \lambda r^2 \\
\\
& \leq & 
\la \nabla f(x) + s, x - T \ra - \lambda r^2 \\
\\
& \leq & 
r \| \nabla f(x) + s \|_{*} - \lambda r^2.
\ea
\eeq
Therefore, we obtain the following inequality.
\BL
For any $s \in \partial \psi(x)$, it holds 
\beq \label{RBound}
\ba{rcl}
r & \leq & \frac{1}{\lambda + \mu} \|\nabla f(x) + s\|_{*}.
\ea
\eeq
\EL
\proof
Indeed, by bounding the Hessian in~\eqref{HessBound1} from below, we get
$$
\ba{rcl}
\mu r^2 & \leq & r \|\nabla f(x) + s\|_{*} - \lambda r^2. \QF
\ea
$$

By maximizing the right hand side of~\eqref{HessBound1} in $r$,
we get the following bound.

\BL
For any $s \in \partial \psi(x)$, it holds 
\beq \label{HessBound2}
\ba{rcl}
\la \nabla^2 f(x)(T - x), T - x \ra & \leq &  
\frac{ \|\nabla f(x) + s\|_{*}^2}{4\lambda}.
\ea
\eeq
\EL

\begin{algorithm}[t]
	\caption{Gradient Regularization of Newton Method}
	\label{alg:FixedNewton}
\begin{algorithmic}[1]
	\Require $x_0 \in \dom \psi$, $\psi'(x_0) \in \partial \psi(x_0)$, $q\in[2, 4]$, $M_q>0$
	\For{$k=0,1,\dotsc$}
		\State $g_{k} = \|\nabla f(x_{k}) + \psi'(x_{k})\|_{*}$
		\State $\lambda_k  =  (6 M_q g_k^{q - 2} )^{\frac{1}{q - 1}}$
		\State $x_{k+1} = \argmin\limits_{x} \left\{
\la \nabla f(x_k), x - x_k \ra + \frac{1}{2} \nabla^2 f(x_k)[x - x_k]^2 + \frac{\lambda_k}{2}\|x - x_k\|^2 + \psi(x)
 \right\}$
		\State $\psi'(x_{k+1})  \Def  - \nabla f(x_k)  - \nabla^2 f(x_k)(x_{k+1} - x_k) - \lambda_k B(x_{k+1} - x_k) $
	\EndFor
\end{algorithmic}
\end{algorithm}

\section{Problem Classes}
\label{SectionProblemClasses}

For the differentiable part of our objective function,
let us introduce some smoothness characteristics.
Namely, let us assume that \textit{either}
its Hessian \textit{or} its third derivative
is H\"older continuous.
For that, let us define the following family of constants:
$$
\ba{rcl}
L_{p, \nu} & \Def & \sup\limits_{\substack{x, y \in \E \\x\not=y}}
\Bigl\{ \,
\frac{\| \nabla^p f(x) - \nabla^p f(y)\| }{\|x - y\|^{\nu}} 
\, \Bigr\},
\ea
$$
where $p = 2$ or $p = 3$ and $\nu \in [0, 1]$.
Since
$$
\ba{rcl}
\ln L_{p, \nu} & = & 
\sup\limits_{\substack{x, y \in \E \\x\not=y}}
\Bigl\{ \,
\ln \| \nabla^p f(x) - \nabla^p f(y)\| 
- \nu \ln \|x - y\|
\, \Bigr\},
\ea
$$
we see that $L_{p, \nu}$ is a \textit{log-convex} function of $\nu$.
Therefore, for any $0 \leq \nu_1 \leq \nu_2 \leq 1$, we have
$$
\ba{rcl}
L_{p, \nu} & \leq & 
\bigl[ L_{p, \nu_1} \bigr]^{\frac{\nu_2 - \nu}{\nu_2 - \nu_1}}
\bigl[ L_{p, \nu_2} \bigr]^{\frac{\nu - \nu_1}{\nu_2 - \nu_1}},
\qquad \forall\nu \in [\nu_1, \nu_2].
\ea
$$
In particular, if $L_{p, 0}$ and $L_{p, 1}$ are finite, we conclude
that all intermediate H\"older constants also exist and 
\beq \label{HolderLogConvex}
\ba{rcl}
L_{p, \nu} & \leq & L_{p, 0}^{1 - \nu} L_{p, 1}^{\nu}, \qquad \forall \nu \in [0, 1].
\ea
\eeq

By the triangle inequality,
$\|\nabla^3 f(x) - \nabla^3 f(y)\| \leq \| \nabla^3 f(x) \| + \| \nabla^3 f(y) \|$.
Thus, for the extreme values of parameters, we have
\beq \label{ValueBorder}
\ba{rcl}
L_{3, 0} & \leq & 2 L_{2, 1}.
\ea
\eeq
Let us consider the following example, which shows that
$L_{2, 1}$ and $L_{3, 0}$ can be different. 
\BE \label{ExampleDiscontinuity}
Let $f(x) = \frac{1}{2}x^2 + \frac{1}{6}|x|^3 \, : \, \R \to \R$. Then,
$$
\ba{rcl}
f'(x) & = & x + \frac{1}{2}|x|x, \quad
f''(x) \; = \;  1 + |x|,
\ea
$$
and for all $x \in \R \setminus \{ 0 \}$, we have $f'''(x)  =  \sign(x)$.
Therefore,
$$
\ba{rcl}
| f''(x) - f''(y) |
& = & \bigl| |x| - |y| \bigr|
\;\; \leq \quad |x - y|
\ea
$$
and
$$
\ba{rcl}
| f'''(x) - f'''(y) | & = & \bigl| \sign(x) - \sign(y) \bigr|
\;\; \leq \;\; 2.
\ea
$$
Hence, $L_{2, 1} = 1$ and $L_{3, 0} = 2$.  \qed
\EE

We have the following useful bound for the right-hand side of \eqref{ValueBorder}.
\BL 
For any $\gamma, \nu \in [0, 1]$, it holds:
\beq \label{BoundL21}
\ba{rcl}
L_{2, 1} & \leq & \frac{2 + \nu}{1 + \nu}
\bigl[L_{2, \gamma} \bigr]^{\frac{\nu}{1 + \nu - \gamma}}
\bigl[ L_{3, \nu}  \bigr]^{\frac{1 - \gamma}{1 + \nu - \gamma}}.
\ea
\eeq
In particular, for $\gamma = 0$ and $\nu = 1$, we get
\beq \label{BoundL01}
\ba{rcl}
L_{2, 1} & \leq & \frac{3}{2} \sqrt{ L_{2, 0} L_{3, 1}}.
\ea
\eeq
\EL
\proof
Indeed, we assume that the third derivative is H\"older continuous of degree $\nu$
with constant $L_{3, \nu} < +\infty$. Hence, 
for all $x, y \in \E$, we have
$$
\ba{rcl}
\| \nabla^2 f(y) - \nabla^2 f(x) - \nabla^3 f(x)[y - x] \| 
& \leq & \frac{L_{3, \nu}}{1 + \nu}\|y - x\|^{1 + \nu}.
\ea
$$
By triangle inequality,
$$
\ba{rcl}
\| \nabla^3 f(x)[y - x] \| & \leq & \| \nabla^2 f(x) - \nabla^2 f(y) \| 
+ \frac{L_{3, \nu}}{1 + \nu}\|y - x\|^{1 + \nu} \\
\\
& \leq & L_{2, \gamma}\|y - x\|^{\gamma}
+ \frac{L_{3, \nu}}{1 + \nu}\|y - x\|^{1 + \nu}.
\ea
$$
Let us take $y := x + \tau h$, where $\|h\| = 1$ and $\tau > 0$. Then,
$$
\ba{rcl}
\| \nabla^3 f(x) \| & \leq & \frac{L_{2, \gamma}}{\tau^{1 - \gamma}}
+ \frac{L_{3, \nu}}{1 + \nu} \tau^{\nu}, \qquad \forall \tau > 0.
\ea
$$
It remains to substitute 
$\tau := \Bigl[\frac{L_{2, \gamma}}{L_{3, \nu}}\Bigr]^{\frac{1}{1 + \nu - \gamma}}$,
which balances both terms.
\qed

We see that parameters of our problem classes
for different $p \in \{2, 3\}$ and $\nu \in [0, 1]$ are related to each other.
It is convenient to have for them a univariate parametrization.
Let us define a family of constants $ M_q \in \R \cup \{ +\infty \}$
with $2 \leq q \leq 4$ as follows:
$$
\ba{rcl}
M_{2 + \nu} & \Def & L_{2, \nu}, \qquad \nu \in [0, 1), \\
\\
M_{3 + \nu} & \Def & L_{3, \nu}, \qquad \nu \in [0, 1].
\ea
$$

Note that by combining~\eqref{HolderLogConvex},\eqref{ValueBorder}, and \eqref{BoundL01},
we obtain
\beq \label{MqBound}
\ba{rcl}
M_q & \leq & 3 
\bigl[ M_2 \bigr]^{\frac{4 - q}{2}} \bigl[ M_4 \bigr]^{\frac{q - 2}{2}}, \qquad \forall q \in [2, 4].
\ea
\eeq
Hence, if $M_2 < +\infty$ and $M_4 < +\infty$, then the whole family is bounded.
Clearly, it can be discontinuous at $q = 3$ (see Example~\ref{ExampleDiscontinuity}).

Our main assumption is as follows.
\BAS
The value $M_q$ is finite for at least one $q\in [2, 4]$:
\beq \label{MainSmooth}
\ba{rcl}
\inf\limits_{2 \leq q \leq 4} M_q & < & +\infty.
\ea
\eeq
\EAS

Note that $M_2 \leq \sup_x \| \nabla^2 f(x) \|$. Thus, we cover
even the standard class for the first-order methods.

It appears that the global complexity of the regularized Newton method
depends on the values $M_q$, $2 \leq q \leq 4$, in a very natural and universal way.
At the same time, it is important that our \textit{super-universal algorithm},
presented in Section~\ref{SectionUniversal}, does not need explicit values of these constants.

\section{Gradient Regularization}
\label{SectionGradReg}

At each step of our method,
we are going to use the following choice of the regularization parameter:
\beq \label{LambdaChoice}
\boxed{
\ba{rcl}
\lambda & := & H g^{\alpha}
\ea
}
\eeq
where $H > 0$ and $\alpha \in [0, 1]$ are some constants, and $g := \|\nabla f(x) + s\|_{*}$
for some $s \in \partial \psi(x)$.
Let us investigate some properties of this choice, 
taking into account our smoothness condition~\eqref{MainSmooth}.
We start with the case $2 \leq q < 3$ (H\"older continuity of the Hessian).

\BT \label{Th1}
Let $g = \| \nabla f(x) + s \|_{*}  > 0$.
Assume that for some $q$, $2 \leq q < 3$,
with $M_q < +\infty$,
our parameters satisfy the following conditions:
\beq \label{AlphaBound1}
\ba{rcl}
\frac{q - 2}{q - 1} & \leq & \alpha \;\; \leq \;\; 1,
\ea
\eeq
and
\beq \label{HBound1}
\ba{rcl}
H & \geq & \Bigl( \frac{1}{q - 1} M_q \Big)^{\frac{1}{q - 1}}
 \Bigl( \frac{1}{g} \Bigr)^{\alpha  - \frac{q - 2}{q - 1}}.
\ea
\eeq
Then, 
\beq \label{Method1}
\ba{rcl}
\la F'(T), x - T \ra & \geq & \frac{1}{2\lambda} \|F'(T)\|_{*}^2.
\ea
\eeq
\ET
\proof
Condition $M_q < +\infty$ implies that the Hessian is H\"older continuous
of degree $\nu = q - 2$ with constant $L_{2, \nu}$.
Thus, we have the following bound for the gradient,
$$
\ba{rcl}
\frac{L_{2, \nu} r^{1 + \nu}}{1 + \nu}
& \geq &
\| \nabla f(T) - \nabla f(x) - \nabla^2 f(x)(T - x)\|_{*} \\
\\
& \overset{\eqref{SubgDef}}{=} &
\| F'(T) + \lambda B(T - x) \|_{*}.
\ea
$$
Squaring both sides of this inequality, we get
$$
\ba{rcl}
\Bigl( \frac{L_{2, \nu} r^{1 + \nu}}{ 1 + \nu} \Bigr)^2
& \geq & 
\| F'(T) \|_{*}^2 + \lambda^2 r^2 + 2\lambda \la F'(T), T - x \ra.
\ea
$$
This means that
$$
\ba{rcl}
\la F'(T), x - T \ra & \geq & \frac{1}{2 \lambda} \|F'(T)\|_{*}^2 
+ \frac{\lambda r^2}{2} - \frac{1}{2\lambda} \Bigl( \frac{L_{2, \nu} r^{1 + \nu}}{ 1 + \nu} \Bigr)^2.
\ea
$$
Hence, for proving~\eqref{Method1}, it is enough to justify
the following relation:
$$
\ba{rcl}
\frac{\lambda r^2}{2} & \geq & \frac{L_{2, \nu}^2 r^{2(1 + \nu)}}{ 2\lambda (1 + \nu)^2}
\quad \Leftrightarrow \quad
\lambda \;\; \geq \;\; \frac{1}{1 + \nu}  L_{2, \nu} r^{\nu}.
\ea
$$
In view of~\eqref{RBound}, 
it is ensured by inequality
$$
\ba{rcl}
\lambda & \geq & 
\frac{1}{1 + \nu} L_{2, \nu} \bigl( \frac{g}{\lambda} \bigr)^{\nu},
\ea
$$
which is equivalent to
$$ 
\ba{rcl}
 \lambda & = & H g^{\alpha}
 \;\; \geq \;\;
\Bigl(   \frac{1}{1 + \nu} L_{2, \nu} g^{\nu} \Bigr)^{\frac{1}{1 + \nu}}
\quad  \Leftrightarrow
\ea
$$
$$
\ba{rcl}
H & \geq &
\bigl(  \frac{1}{1 + \nu} L_{2, \nu} \bigr)^{\frac{1}{1 + \nu} }
\cdot \Bigl( \frac{1}{g}  \Bigr)^{\alpha - \frac{\nu}{1 + \nu}}
\;\; = \;\;
\bigl( \frac{1}{q - 1} M_q  \bigr)^{\frac{1}{q - 1}}
\cdot \Bigl( \frac{1}{g} \Bigr)^{\alpha - \frac{q - 2}{q - 1}}. \QF
\ea
$$

Now, let us analyze the case $3 \leq q \leq 4$.
This is H\"older continuity of the third derivative
with parameter $\nu = q - 3$.
Firstly, let us bound the third derivative of $f(\cdot)$ along
direction $T - x$.

\BL \label{LemmaThird}
 For any $\nu \in [0, 1]$ and $s \in \partial  \psi(x)$, we have
\beq \label{ThirdDNormBound}
\ba{rcl}
\| \nabla^3 f(x)[T - x]^2 \|_{*} & \leq & 
2 \Bigl(  \frac{1}{1 + \nu} L_{3, \nu} r^{2} \Bigr)^{\frac{1}{1 + \nu}}
\Bigl(  \frac{\| \nabla f(x) + s  \|_{*}^2}{4\lambda} \Bigr)^{\frac{\nu}{1 + \nu}}.
\ea
\eeq
\EL
\proof
By convexity of $f(\cdot)$ and H\"older continuity of the third derivative, 
for all $x, y, h \in \E$, we have
$$
\ba{rcl}
0 & \leq & \nabla^2 f(y)[h]^2
\;\; \leq \;\; \nabla^2 f(x)[h]^2 
+ \la \nabla^3 f(x)[h]^2, y - x\ra
+ \frac{L_{3, \nu} \|y - x\|^{1 + \nu}}{1 + \nu} \|h\|^2.
\ea
$$
Taking $y = x + \tau u$,
for $u \in \E$ with $\|u\| = 1$ and arbitrary $\tau > 0$, we get
$$
\ba{rcl}
\| \nabla^3 f(x)[h]^2 \|_{*}
& \leq & \frac{1}{\tau} \nabla^2 f(x)[h]^2
+ \frac{\tau^{\nu} L_{3, \nu}}{1 + \nu} \|h\|^2.
\ea
$$
For $\tau := \Bigl(
\frac{ (1 + \nu) \nabla^2 f(x)[h]^2}{ L_{3, \nu} \|h\|^2 }
\Bigr)^{\frac{1}{1 + \nu}}$, this gives
$$
\ba{rcl}
\| \nabla^3 f(x)[h]^2 \|_{*}
& \leq & 
2 \Bigl(  \frac{1}{1 + \nu} L_{3, \nu} \|h\|^2  \Bigr)^{\frac{1}{1 + \nu}}
\Bigl(  \nabla^2 f(x)[h]^2 \Bigr)^{\frac{\nu}{1 + \nu}}.
\ea
$$
Choosing $h := T - x$ , we get the desired 
inequality by \eqref{HessBound2}. 
\qed

Let us prove now a lower bound for the progress in one iteration.

\BT \label{Th2}
Let $g = \| \nabla f(x) + s \|_{*}  > 0$.
Assume that for some $q$, $3 \leq q \leq 4$,
with $M_q < +\infty$,
our parameters satisfy the following conditions:
\beq \label{AlphaBound2}
\ba{rcl}
\frac{q - 2}{q - 1} & \leq & \alpha \;\; \leq \;\; 1,
\ea
\eeq
and
\beq \label{HBound2}
\ba{rcl}
H & \geq & \Bigl( \frac{6^{q - 2}}{4^{q - 3}(q - 2)}  M_q \Big)^{\frac{1}{q - 1}}
\Bigl( \frac{1}{g} \Bigr)^{\alpha - \frac{q - 2}{q - 1}}.
\ea
\eeq
Then, 
\beq \label{Method2}
\ba{rcl}
\la F'(T), x - T \ra & \geq & \frac{1}{4\lambda} \|F'(T)\|_{*}^2.
\ea
\eeq
\ET
\proof
Condition $M_q < +\infty$ implies that the third derivative is
H\"older continuous of degree $\nu = q - 3$ with constant $L_{3, \nu}$.
Hence,
$$
\ba{rcl}
\frac{L_{3, \nu} r^{2 + \nu}}{(1 + \nu)(2 + \nu)}
& \geq & 
\| \nabla f(T) - \nabla f(x) - \nabla^2 f(x)(T - x) - \frac{1}{2}\nabla^3 f(x)[T - x]^2 \|_{*} \\
\\
& \overset{\eqref{SubgDef}}{=} &
\| F'(T) + \lambda B(T - x) - \frac{1}{2}\nabla^3 f(x)[T - x]^2 \|_{*}.
\ea
$$
Squaring both sides of this inequality, we get
$$
\ba{rcl}
\Bigl(  \frac{L_{3, \nu} r^{2 + \nu}}{(1 + \nu)(2 + \nu)} \Bigr)^2
& \geq &
\| F'(T) \|_{*}^2 + \lambda^2 r^2 + \frac{1}{4} \| \nabla^3 f(x)[T - x]^2  \|_{*}^2  \\
\\
& & \; +  \;2\lambda \la F'(T), T - x \ra
-  \la \nabla^3 f(x)[T - x]^2, B^{-1} F'(T) \ra \\
\\
& & \; - \; \lambda \la \nabla^3 f(x)[T - x]^2, T - x \ra \\
\\
& \geq &
\frac{1}{2}\| F'(T) \|_{*}^2 + \lambda^2 r^2  - \frac{1}{4}\| \nabla^3 f(x) [T - x]^2 \|_{*}^2  \\
\\
& & \; + \;
 2\lambda \la F'(T), T - x \ra - \lambda r \| \nabla^3 f(x)[T - x]^2 \|_{*},
\ea
$$
where we applied Cauchy-Schwartz and Young's inequalities in the last bound. 
Then, rearranging the terms and using Lemma~\ref{LemmaThird}, we have
$$
\ba{rcl}
\la F'(T), x - T \ra & \geq & 
\frac{1}{4\lambda} \|F'(T)\|_{*}^2
+ \frac{\lambda r^2}{2}
- \frac{1}{8\lambda}\| \nabla^3 f(x)[T - x]^2\|_*^2 \\
\\
& & 
\; - \; \frac{r}{2} \| \nabla^3 f(x)[T - x]^2\|_{*}
- \frac{1}{2\lambda}\Bigl(  \frac{L_{3, \nu} r^{2 + \nu}}{(1 + \nu)(2 + \nu)} \Bigr)^2 \\
\\
& \overset{\eqref{ThirdDNormBound}}{\geq} &
\frac{1}{4\lambda} \|F'(T)\|_{*}^2
+ \frac{\lambda r^2}{2}
- \frac{1}{2\lambda}\Bigl( \frac{L_{3, \nu} r^2}{1 + \nu}  \Bigr)^{\frac{2}{1 + \nu}}
\Bigl(  \frac{g^2}{4\lambda}  \Bigr)^{\frac{2\nu}{1 + \nu}} \\
\\
& &
\; - \;
r  \Bigl(
\frac{L_{3, \nu} r^{2} }{1 + \nu} 
\Bigr)^{\frac{1}{1 + \nu}}
\Bigl( \frac{g^2}{4 \lambda} \Bigr)^{\frac{\nu}{1 + \nu}}
- \frac{1}{2\lambda}\Bigl(  \frac{L_{3, \nu} r^{2 + \nu}}{(1 + \nu)(2 + \nu)} \Bigr)^2.
\ea
$$
Let us divide the term $\frac{\lambda r^2}{2}$
into three equal parts.
Then we need to ensure validity of three inequalities.
\begin{enumerate}
	\item $\frac{\lambda r^2}{6} \geq 
	\frac{1}{2\lambda} \Bigl( \frac{L_{3, \nu} r^2}{1 + \nu}  \Bigr)^{\frac{2}{1 + \nu}}
	\Bigl(  \frac{g^2}{4\lambda}  \Bigr)^{\frac{2\nu}{1 + \nu}}$
	$\quad \Leftrightarrow \quad $ 
	$
	\lambda \geq \frac{3}{\lambda} 
	\Bigl(  \frac{L_{3, \nu}}{1 + \nu} \Bigr)^{\frac{2}{1 + \nu}}
	\Bigl(  \frac{g^2}{4\lambda} \Bigr)^{\frac{2\nu}{1 + \nu}}
	r^{\frac{2(1 - \nu)}{1 + \nu}}.
	$

	In view of~\eqref{RBound}, a sufficient condition is
	$$
	\ba{rcl}
	\lambda & \geq & 
	\frac{3}{\lambda} 
	\Bigl(  \frac{L_{3, \nu}}{1 + \nu} \Bigr)^{\frac{2}{1 + \nu}}
	\Bigl(  \frac{g^2}{4\lambda} \Bigr)^{\frac{2\nu}{1 + \nu}}
	\Bigl( \frac{g}{\lambda} \Bigr)^{\frac{2(1 - \nu)}{1 + \nu}},
	\ea
	$$
	which is equivalent to
	$
	H  \; = \; \frac{\lambda}{g^{\alpha}}
	\; \geq \;
	\Bigl(
	\frac{3^{(1 + \nu) / 2} L_{3, \nu}}{4^{\nu} (1 + \nu)} 
	\Bigr)^{\frac{1}{2 + \nu}}
	\Bigl( \frac{1}{g} \Bigr)^{\alpha - \frac{1 + \nu}{2 + \nu}}.
	$
	
	\item $\frac{\lambda r^2}{6} \geq 
	r  \Bigl(
	\frac{L_{3, \nu} r^{2} }{1 + \nu} 
	\Bigr)^{\frac{1}{1 + \nu}}
	\Bigl( \frac{g^2}{4 \lambda} \Bigr)^{\frac{\nu}{1 + \nu}}$
	$\quad \Leftrightarrow \quad$ 
	$
	\lambda \geq 6 \Bigl( \frac{L_{3, \nu}}{(1 + \nu)} \Bigr)^{\frac{1}{1 + \nu}}
	\Bigl( \frac{g^2}{4 \lambda}  \Bigr)^{\frac{\nu}{1 + \nu}} r^{\frac{1 - \nu}{1 + \nu}}.
	$
	
	In view of~\eqref{RBound}, a sufficient condition is
	$$
	\ba{rcl}
	\lambda & \geq & 6 \Bigl( \frac{L_{3, \nu}}{(1 + \nu)} \Bigr)^{\frac{1}{1 + \nu}}
	\Bigl( \frac{g^2}{4 \lambda}  \Bigr)^{\frac{\nu}{1 + \nu}} 
	\Bigl(  \frac{g}{\lambda} \Bigr)^{\frac{1 - \nu}{1 + \nu}},
	\ea
	$$
	which is equivalent to
	$
	H \; = \; \frac{\lambda}{g^{\alpha}}
	\; \geq \;
	\Bigl( 
	\frac{6^{1 + \nu} L_{3, \nu}}{4^{\nu}(1 + \nu)}
	\Bigr)^{\frac{1}{2 + \nu}} 
	\Bigl( \frac{1}{g} \Bigr)^{\alpha - \frac{1 + \nu}{2 + \nu}}.
	$
	
	\item
	$ \lambda \geq 
	\frac{3}{\lambda} \Bigl( \frac{L_{3, \nu}}{(1 + \nu)(2 + \nu)} \Bigr)^2 r^{2(1 + \nu)}$.
	Hence, due to~\eqref{RBound}, a sufficient condition is
	$$
	\ba{rcl}
	\lambda & \geq & 
	\frac{3}{\lambda} \Bigl(  \frac{L_{3, \nu}}{(1 + \nu)(2 + \nu)} \Bigr)^2
	\Bigl( \frac{g}{\lambda} \Bigr)^{2(1 + \nu)},
	\ea
	$$
	which is equivalent to
	$
	H \; = \; \frac{\lambda}{g^{\alpha}}
	\; \geq \;
	\Bigl(  \frac{3^{1/2} L_{3, \nu}}{(1 + \nu)(2 + \nu)} \Bigr)^{\frac{1}{2 + \nu}}
	\Bigl(  \frac{1}{g} \Big)^{\alpha - \frac{1 + \nu}{2 + \nu}}.
	$
\end{enumerate}
We see that in all three cases, the lower bounds for $H$ are very similar.
Thus, it is sufficient to choose one with the maximal absolute constant.
This is 
$$
\ba{rcl}
H 
& \geq &
\Bigl(
\frac{6^{1 + \nu} L_{3, \nu}}{4^{\nu} (1 + \nu)} 
\Bigr)^{\frac{1}{2 + \nu}}
\Bigl( \frac{1}{g} \Bigr)^{\alpha - \frac{1 + \nu}{2 + \nu}}
\;\; = \;\;
\Bigl( 
\frac{6^{q - 2}}{4^{q - 3}(q - 2)} M_q
\Bigr)^{\frac{1}{q - 1}} 
\Bigl( \frac{1}{g} \Bigr)^{\alpha - \frac{q - 2}{q - 1}}. \QF
\ea
$$

\BC \label{MainCor}
Let $g = \| \nabla f(x) + s \|_{*} > 0$
and $M_q < +\infty$ for some $q \in [2, 4]$.
Then, for
$
\frac{q - 2}{q - 1} \leq  \alpha \leq 1,
$
and
$
H \geq \bigl( 6 M_q \bigr)^{\frac{1}{q - 1}} 
\Bigl( \frac{1}{g} \Bigr)^{\alpha - \frac{q - 2}{q - 1}},
$
it holds
\beq \label{MainProgress}
\ba{rcl}
\la F'(T), x - T \ra & \geq & \frac{1}{4\lambda} \|F'(T)\|_{*}^2.
\ea
\eeq
\EC
\proof
Indeed, for $2 \leq q \leq 3$, we have $\frac{1}{q - 1} \leq 6$,
and for $3 \leq q \leq 4$, we also have
$\frac{6^{q - 2}}{4^{q - 3} (q - 2)} \leq 6$.
\qed

\BR
Note that inequality \eqref{MainProgress} implies
\beq \label{NewGradBound}
\ba{rcl}
g_{+} \;\; \Def \;\; \| F'(T) \|_{*}
& \leq & 4 \lambda r
\;\; \overset{\eqref{RBound}}{\leq}
4g.
\ea
\eeq
\ER

Now, let us look at the simplest way of choosing regularization constants,
when parameter $q \in [2, 4]$ is known and fixed.
By Corollary~\ref{MainCor}, we can take
$$
\ba{rcl}
\alpha & := & \frac{q - 2}{q - 1} 
\qquad \text{and} \qquad
H \;\; := \;\; 
\bigl( 6 M_q \bigr)^{\frac{1}{q - 1}}.
\ea
$$
This way, we obtain Algorithm~\ref{alg:FixedNewton}.

By convexity,
we get the following progress for one step of this method:
\beq \label{BasicOneStep}
\ba{rcl}
F(x_k) - F(x_{k + 1}) & \geq &
\la F'(x_{k + 1}), x_k - x_{k + 1} \ra \\
\\
& \overset{\eqref{MainProgress}}{\geq} &
\frac{g_{k + 1}^2}{4 \lambda_k}
\;\; = \;\;
\frac{1}{4 (6M_q)^{1 / (q - 1) }}
\bigl( \frac{g_{k + 1}}{g_k} \bigr)^2
g_k^{\frac{q}{q - 1}}.
\ea
\eeq
This inequality results in a \textit{global convergence rate}
for our process. 
In the next Section~\ref{SectionUniversal},
we derive it explicitly.
However, the main drawback of this scheme
is that we need to fix the degree of smoothness $q \in [2, 4]$
in advance. The parameter $M_q$ is also needed.
Hence, the above scheme is completely theoretical and
cannot be used in practice.
The \textit{super-universal} method,
presented in Section~\ref{SectionUniversal},
resolves both these issues by a simple search procedure.

\section{Super-Universal Method}
\label{SectionUniversal}

At each iteration of this scheme, we adjust the regularization
constant $H_k$ for ensuring inequality \eqref{MainProgress}.
Degree $\alpha \in [\frac{2}{3}, 1]$ of the gradient regularization 
is chosen in advance and does not depend on
a particular problem class.

\begin{algorithm}[t]
	\caption{Super-Universal Newton Method}
	\label{alg:AdaptiveNewton}
\begin{algorithmic}[1]
	\Require $x_0 \in \dom \psi$, $\psi'(x_0) \in \partial \psi(x_0)$.
	Choose arbitrary $\alpha\in\bigl[\frac{2}{3}, 1\bigr]$, $H_0>0$
	\For{$k=0,1,\dotsc$}
		\State $g_{k} = \|\nabla f(x_{k}) + \psi'(x_k)\|_{*}$
		\On{$j_k=0,1,\dotsc$}
			\State $\lambda_k = 4^{j_k} H_k g_k^{\alpha}$
			\State $x_+ = \argmin\limits_{x} \left\{
\la \nabla f(x_k), x - x_k \ra + \frac{1}{2} \nabla^2 f(x_k)[x - x_k]^2 + \frac{\lambda_k}{2}\|x - x_k\|^2 + \psi(x)
 \right\}$
 			\State $\psi'(x_{+}) \Def  - \nabla f(x_k) - \nabla^2 f(x_k)(x_{+} - x_k) - \lambda_k B(x_{+} - x_k)$ 
 			\State $F'(x_{+}) \Def \nabla f(x_{+}) + \psi'(x_{+})$
		\EndOn{$\la F'(x_{+}), x_k - x_{+} \ra \geq  \frac{\| F'(x_{+})\|^2_*}{4 \lambda_k}$}
		\State $x_{k+1} = x_+$
		\State $H_{k + 1} = \frac{4^{j_k} H_k}{4}$
	\EndFor
\end{algorithmic}
\end{algorithm}

We need to prove first that the method is well-defined. Denote
\beq \label{DefH}
\ba{rcl}
\mathcal{H}_{\alpha}(t)
& \Def &
\inf\limits_{2 \leq q \leq 4} \,
 \bigl( 6 M_q \bigr)^{\frac{1}{q - 1}} 
\Bigl( \frac{1}{t} \Bigr)^{\alpha - \frac{q - 2}{q - 1}}, \qquad t > 0.
\ea 
\eeq
Since $2 \leq q \leq 4$ and $\alpha \geq \frac{2}{3}$,
this function is decreasing in $t$.

\BL
Assume that $M_q < +\infty$ for some $q \in [2, 4]$
and
\beq \label{H_0_Bound}
\ba{rcl}
H_0 & \leq & \mathcal{H}_{\alpha}(g_0).
\ea
\eeq
Let for all 
iterations $\{ x_i \}_{i = 0}^{k - 1}$
of Algorithm~\ref{alg:AdaptiveNewton} with some $k \geq 1$, we have
$$
\ba{rcl}
g_i \;\; \Def \;\; \|F'(x_i)\|_{*} & > & 0.
\ea
$$
Then, 
\beq \label{H_k_Bound}
\ba{rcl}
H_{i + 1} & \leq &
\mathcal{H}_{\alpha}(g_{i}),
\qquad 0 \leq i \leq k-1.
\ea
\eeq
Moreover, the total number $N_k$ of oracle calls during the first $k$ iterations
is bounded as follows:
\beq \label{TotalItersBound}
\ba{rcl}
N_k & \leq & 2k +  \frac{1}{2} \log_2 \frac{\mathcal{H}_{\alpha}(g_{k - 1})}{H_0}.
\ea
\eeq
\EL
\proof
Let us prove~\eqref{H_k_Bound} by induction.
Denote formally $g_{-1} \Def g_0$, and then \eqref{H_0_Bound} 
is the base of the induction.
Now, consider the $i$th iteration of the method.

In case $j_i > 0$, the condition of the search procedure is not satisfied
for the previous $\lambda := 4^{j_i - 1}H_i g_i^{\alpha} = H_{i + 1} g_i^{\alpha}$.
Hence, by Corollary~\ref{MainCor}, we conclude that
$$
\ba{rcl}
H_{i + 1} & \leq & \mathcal{H_{\alpha}}(g_{i}).
\ea
$$ 
In the other case, we have $j_i = 0$, and 
$$
\ba{rcl}
H_{i + 1} & = & \frac{H_i}{4}
\;\; \leq \;\; 
\frac{1}{4}\mathcal{H_{\alpha}}(g_{i - 1})
\;\; \overset{\eqref{NewGradBound}}{\leq} \;\;
\frac{4^{\alpha - \frac{q - 2}{q - 1}}}{4}
\mathcal{H_{\alpha}}(g_{i})
\;\; \leq \;\;
\mathcal{H_{\alpha}}(g_{i}).
\ea
$$
Thus, \eqref{H_k_Bound} holds for all $0 \leq i \leq k - 1$.

In order to estimate the total number of oracle calls,  note that $4H_{k + 1} = 4^{j_k} H_k$, where $j_k$ is the number of unsuccessful inner-loop iterations of Algorithm~\ref{alg:AdaptiveNewton}.
Hence, we have
$$
\ba{rcl}
N_k & = & 
\sum\limits_{i = 0}^{k - 1} ( 1  + j_i )
\;\; = \;\;
k + \sum\limits_{i = 0}^{k - 1}  \log_4 \frac{4 H_{i + 1}}{H_i}
\;\; = \;\; 
2k + \log_4 H_k - \log_4 H_0  \\
\\
& \overset{\eqref{H_k_Bound}}{\leq} &
2k + \log_4 \frac{\mathcal{H}_{\alpha}(g_{k - 1})}{H_0}
\;\; = \;\;
2k + \frac{1}{2} \log_2 \frac{\mathcal{H}_{\alpha}(g_{k - 1})}{H_0}.
\QR
\ea
$$

Substituting the bound~\eqref{H_k_Bound} into the formula 
for our choice of $\lambda_k$, we get
\beq \label{Lambda_k_Bound}
\ba{rcl}
\lambda_k & = & 4 H_{k + 1} g_k^{\alpha}
\;\; \overset{\eqref{H_k_Bound}}{\leq} \;\;
4 \mathcal{H}_{\alpha}(g_k) g_k^{\alpha} 
\;\; \overset{\eqref{DefH}}{\leq} \;\;
4 \bigl( 6 M_q \bigr)^{\frac{1}{q - 1}} g_k^{\frac{q - 2}{q - 1}},
\ea
\eeq
for any $q \in [2, 4]$.
Thus, one iteration of our adaptive scheme
ensures that
\beq \label{AdaptiveOneStep}
\ba{rcl}
F(x_k) - F(x_{k + 1})
& \geq & 
\la F'(x_{k + 1}), x_k - x_{k + 1} \ra
\;\; \geq \;\; \frac{g_{k + 1}^2}{4 \lambda_k} \\
\\
& \overset{\eqref{Lambda_k_Bound}}{\geq} &
\frac{1}{16 (6 M_q)^{1 / (q - 1)}}
\bigl( \frac{g_{k + 1}}{g_k} \bigr)^2
g_{k}^{\frac{q}{q - 1}}.
\ea
\eeq
Up to the factor $\frac{1}{4}$,
this bound is the same as inequality \eqref{BasicOneStep}
for the basic method.
However, our new method 
is \textit{adaptive} and it does not need to know any particular values
of $q$ and $M_q$.

Note that the parameter $\alpha$ 
can be chosen \underline{\textit{arbitrarily}}
in the interval $[\frac{2}{3}, 1]$.
For example, one can stick to the choice $\alpha = 1$.
As we see from~\eqref{TotalItersBound}, the price for the universality
is, on average, just \textit{one extra oracle call} per iteration.

By the initial condition \eqref{H_0_Bound}, $H_0$ has to be small.
In fact, this requirement is not restrictive. 
For fulfilling it, we can start with an arbitrary value for $H_0$ 
and decrease it twice until stopping condition from the search procedure is satisfied.
There are two options: either the condition is violated at some moment, and
hence $H_0$ satisfies \eqref{H_0_Bound} by Corollary~\ref{MainCor},
or the gradient becomes smaller and smaller with a linear rate. 
Thus this simple search is of the logarithmic length,
and it can be used at a preliminary stage.

We are ready to prove the global rate of convergence of Algorithm~\ref{alg:AdaptiveNewton}.
Denote by 
\[
\ba{rcl}
\mathcal{F}_0  \; \Def \;  
\Bigl\{ x  \in \dom \psi \; : \;\;  F(x) \leq F(x_0)  \Bigr\}
\ea
\]
the initial sublevel set,
which we assume to be bounded:
$$
\ba{rcl}
D & \Def & \sup\limits_{x, y \in \mathcal{F}_0} 
\|x - y\|
 \;\; < \;\; +\infty.
\ea
$$
By convexity of $F(\cdot)$ and monotonicity of the sequence $\{ F(x_k) \}_{k \geq 0}$,
we have
\beq \label{ConvexityFD}
\ba{rcl}
g_k & \geq & \frac{F_k}{D},
\quad 
F_k \;\; \Def \;\; F(x_k) - F_*.
\ea
\eeq
Without loss of generality, we can assume $F_k > 0$ for all $k \geq 0$.

For $0 \leq \beta \leq 1$, function $y(x) = x^{\beta}, x \geq 0$
is concave, which implies
\beq \label{Concavity1}
\ba{rcl}
a^{\beta} - b^{\beta} & \geq & 
\frac{\beta}{a^{1 - \beta}} (a - b), \qquad \forall a > b \geq 0.
\ea
\eeq
Thus, for $\beta := \frac{1}{q - 1} \in [\frac{1}{3}, 1]$, we have
\beq \label{ConvexOneStep}
\ba{cl}
& \frac{1}{F_{k + 1}^{\beta}}
- \frac{1}{F_k^{\beta}}
\;\; = \;\;
\frac{F_k^{\beta} - F_{k + 1}^{\beta}}{ F_k^{\beta} F_{k+1}^{\beta} }
\;\; \overset{\eqref{Concavity1}}{\geq} \;\;
\frac{\beta(F_k - F_{k + 1})}{F_k F_{k + 1}^{\beta} } \\
\\
& \overset{\eqref{AdaptiveOneStep}}{\geq}
\frac{\beta}{16(6M_q)^{\beta}}
\bigl( \frac{g_{k + 1}}{g_k}  \bigr)^{2}
\frac{g_{k}^{1 + \beta}}{F_k F_{k + 1}^{\beta} } 
\;\; \overset{\eqref{ConvexityFD}}{\geq} \;\;
\frac{\beta}{16(6M_q)^{\beta} D^{1 + \beta}}
\bigl( \frac{g_{k + 1}}{g_k}  \bigr)^{2}.
\ea
\eeq
Telescoping the last bound and
using the inequality for arithmetic and geometric means, we get
\beq \label{ConvexTelescopedStep}
\ba{cl}
& \frac{1}{F_k^{\beta}} 
-  \frac{1}{F_0^{\beta}}
\;\; \geq \;\;
\frac{\beta}{16(6M_q)^{\beta} D^{1 + \beta}}
\sum\limits_{i = 0}^{k - 1} \bigl( \frac{g_{i + 1}}{g_i} \bigr)^{2} \\
\\
& \; \geq \;
\frac{\beta k}{16(6M_q)^{\beta} D^{1 + \beta}}
\Bigl( \, \prod\limits_{i = 0}^{k - 1} \frac{g_{i + 1}}{g_i} \, \Bigr)^{\!\frac{2}{k}}
\;\; = \;\;
\frac{\beta k}{16(6M_q)^{\beta} D^{1 + \beta}} \bigl( \frac{g_k}{g_0} \bigr)^{\frac{2}{k}}.
\ea
\eeq

These observations prove the following global rate.

\BT \label{TheoremGlobalConvex}
Let $M_q < +\infty$ for some $q \in [2, 4]$
and the initial value $H_0$ satisfies~\eqref{H_0_Bound}.
Then, for all $k \geq 1$, we have
\beq \label{MainConvexRate}
\ba{rcl}
F(x_k) - F_* & \leq &
6 M_q D^q
\Bigl(\frac{ 32(q - 1) }{k}\Bigr)^{q - 1}
+ g_0 D \exp\Bigl( -\frac{k}{4} \Bigr). 
\ea
\eeq
\ET
\proof
Indeed, from~\eqref{ConvexTelescopedStep}, we have
\beq \label{ConvProofTelescoped}
\ba{cl}
& \frac{1}{F_k^{\beta}} - \frac{1}{F_0^{\beta}}
\;\; \overset{\eqref{ConvexTelescopedStep}}{\geq} \;\; 
\frac{\beta k}{16 (6 M_q)^{\beta} D^{1 + \beta}}
\bigl( \frac{g_k}{g_0}   \bigr)^{ \frac{2}{k} }
\overset{\eqref{ConvexityFD}}{\geq} 
\frac{\beta k}{16 (6 M_q)^{\beta} D^{1 + \beta}}
\bigl( \frac{F_k}{g_0 D}   \bigr)^{\frac{2}{k}} \\
\\
& = 
\frac{\beta k}{16 (6 M_q)^{\beta} D^{1 + \beta}}
 \exp\Bigl( 
- \frac{2}{k} \ln \frac{g_0 D}{F_k} 
\Bigr)
 \geq 
\frac{\beta k}{16 (6 M_q)^{\beta} D^{1 + \beta}}
 \Bigl( 1 - \frac{2}{k}  \ln \frac{g_0 D}{F_k} \Bigr).
\ea
\eeq
It remains to consider two cases.
Either
$$
\ba{rcl}
\frac{2}{k} \ln \frac{g_0 D}{F_k}
& \geq & \frac{1}{2}
\quad \Leftrightarrow \quad
F_k \;\; \leq \;\; g_0 D \exp\bigl( -\frac{k}{4} \bigr),
\ea
$$
or $\frac{2}{k} \ln \frac{g_0 D}{F_k} < \frac{1}{2}$,
which together with \eqref{ConvProofTelescoped} leads to
$$
\ba{rcl}
\frac{1}{F_k^{\beta}} & \geq & 
\frac{\beta k}{32 (6 M_q)^{\beta} D^{1 + \beta}}
\quad \Leftrightarrow \quad
F_k \;\; \leq \;\;
\Bigl[\frac{32(6 M_q)^{\beta} D^{1 + \beta}}{\beta k} \Bigr]^{\frac{1}{\beta}}
\;\; = \;\;
\frac{ (32(q - 1))^{q - 1}6 M_q D^{q} }{k^{q - 1}}.
\ea
$$
Combining these two bounds, we get inequality~\eqref{MainConvexRate}.
\qed

Note that the second term in \eqref{MainConvexRate} decreases exponentially in $k$.
Indeed, for any $\varepsilon > 0$, starting from the moment
\beq \label{KLogBound}
\ba{rcl}
k & \geq & 4\ln\frac{g_0 D}{\varepsilon},
\ea
\eeq
this term is bounded by $\varepsilon$.

\BC \label{SummaryCor}
	Assume that $k\ge 4\ln\frac{g_0 D}{\varepsilon}$. If the third derivative is Lipschitz continuous ($q=4$), we obtain the same convergence rate as that of the third-order Tensor Method:
	\[
	\ba{rcl}
		F(x_k) - F_* & \leq & O\bigl( \frac{M_4 D^4}{k^3} \bigr).
	\ea
	\]
	If the Hessian is Lipschitz continuous ($q=3$), then our method achieves the same convergence rate as Cubic Newton Method~\cite{nesterov2006cubic}:
	\[
	\ba{rcl}
		F(x_k) - F_*   & \leq &  O\bigl( \frac{M_3 D^3}{k^2} \bigr).
	\ea
	\]
	Finally, if the Hessian has bounded variation ($q=2$), then the rate is:
	\[
	\ba{rcl}
		F(x_k) - F_*  & \leq &  O\bigl( \frac{M_2 D^2}{k} \bigr).
	\ea
	\]
\EC

The last rate in Corollary~\ref{SummaryCor} is typical for the Gradient Methods~\cite{nesterov2018lectures}. However,
now the constant $M_2$ bounds the \textit{variation} of the Hessian. It can be
much smaller than the norm of the Hessian,
which is used in the analysis of the first-order methods.

As we have seen, depending on the problem class,
the global convergence rate can vary significantly.
Fortunately, our super-universal method does not fix any particular
$q$ and thus achieves the best complexity among all variants.

\BC
According to Theorem~\ref{TheoremGlobalConvex},
in order to reach $F(x_k) - F_* \leq \varepsilon$, it is enough to perform
\beq \label{ConvexTotalComplexity}
\ba{rcl}
k & = & 
O\Bigl( \,
\inf\limits_{2 \leq q \leq 4}
\Bigl[
\frac{M_q D^{q}}{\varepsilon} \Bigr]^{\frac{1}{q - 1}}
+ \ln\frac{g_0 D}{\varepsilon}
\,
\Bigr)
\ea
\eeq
iterations of Algorithm~\ref{alg:AdaptiveNewton}.
\EC

\section{Strictly Convex Functions}
\label{SectionSConvex}

Let us analyze convergence of Algorithm~\ref{alg:AdaptiveNewton}
on some subclasses of \textit{strictly convex} functions.
As we will see, an automatic acceleration on such functions
ensures much faster global rates, as well as the local superlinear convergence.

As in the previous sections, 
consider the initial sublevel set
$$
\ba{rcl}
\mathcal{F}_0 & = & \Bigl\{ x \in \dom\psi : \; F(x) \leq F(x_0) \Bigr\}.
\ea
$$
We have used the primal norm $\| \cdot \|$ 
for measuring its size, denoting
$$
\ba{rcl}
D & \Def & \sup\limits_{x, y \in \mathcal{F}_0} \|x - y\|.
\ea
$$

However, this is not the only possibility.
The other natural measure would be
the \textit{symmetrized Bregman divergence} induced by our objective:
$$
\ba{rcl}
\beta_F(x, y) & \Def & 
\la G_x - G_y, x - y \ra,
\ea
$$ 
for some fixed selection of subgradients 
$G_x \in \partial F(x)$ and $G_y \in \partial F(y)$.
Note that strict convexity ensures $\beta_F(x, y) > 0$ for $x \not= y$.
Defining
\beq \label{VFDef}
\ba{rcl}
V_F & \Def & \sup\limits_{x, y \in \mathcal{F}_0}
\beta_{F}(x, y),
\ea
\eeq
and assuming its boundedness,
we can use the following normalized measure:
$$
\ba{rcl}
\xi_F(x, y) & \Def & 
\frac{1}{V_F} \beta_F(x, y)
\;\; \leq \;\; 1, \qquad \forall x, y \in \mathcal{F}_0.
\ea
$$

It is interesting that the relations between these two measures
have important consequences for complexity of the corresponding
problem~\eqref{MainProblem}.
Let us introduce a new characteristic called
$s$\textit{-relative size} ($s \geq 2$). 
Denote\footnote{We use a conventional notation $1/\infty = 0$.}
\beq \label{RelSize}
\ba{rcl}
D_s & \Def &
\sup\limits_{\substack{x, y \in \mathcal{F}_0 \\ x\not= y}}
\Bigl\{ \;
\|x - y\| \cdot \xi_F(x, y)^{-1/s}
\; \Bigr\}.
\ea
\eeq
Thus, by our definition $D_{\infty} = D$.
Note that
$$
\ba{rcl}
\ln D_s & = &
\sup\limits_{\substack{x, y \in \mathcal{F}_0 \\ x\not= y}} \Bigl\{   \;
\ln \|x - y\|
+ \frac{1}{s} \ln\frac{1}{\xi_F(x, y)}
\; \Bigr\} \\
\\
& = & 
\sup\limits_{\substack{x, y \in \mathcal{F}_0 \\ x\not= y}} \Bigl\{   \;
\ln \|x - y\|
+ \frac{1}{s} \ln\frac{V_F}{\beta_F(x, y)}
\; \Bigr\}.
\ea
$$
Since the last expression is a pointwise supremum of convex functions,
we conclude that $D_s$ is a \textit{log-convex} function of $s$.
Hence, if for some $2 \leq s_1 \leq s_2$
we have $D_{s_i} < +\infty, i \in \{1, 2\}$,
then 
\beq \label{RhoSegmentBound}
\ba{rcl}
D_s & \leq & 
\bigl[ D_{s_1}  \bigr]^{\frac{s_2 - s}{s_2 - s_1}}
\bigl[ D_{s_2} \bigr]^{\frac{s - s_1}{s_2 - s_1}},
\qquad s_1 \leq s \leq s_2,
\ea
\eeq
and $D_s$ is continuous on this segment.

\BE
Let $F(x) = \frac{1}{2}\|x\|^2$. Then $\beta_F(x, y) = \|x - y\|^2$
and $V_F = D^2$. Consequently,
$$
\ba{rcl}
D_s & = &
\sup\limits_{\substack{x, y \in \mathcal{F}_0 \\ x \not= y}}
\Bigl\{
\;
\|x - y\|^{1 - \frac{2}{s}} D^{\frac{2}{s}}
\;
\Bigr\}
\;\; = \;\; D, \qquad  2 \leq s \leq \infty.
\ea
$$
\EE

\BE
Let $F$ be uniformly convex of degree $s \geq 2$. Then
\beq \label{UConvex}
\ba{rcl}
\la G_x - G_y, x - y \ra & \geq & \sigma_s \|x - y\|^s,
\qquad \forall x, y \in \dom F,
\ea
\eeq
for all $G_x \in \partial F(x), G_y \in \partial F(y)$
and some $\sigma_s > 0$. Hence,
$$
\ba{rcl}
D_s
& \leq & \bigl(  \frac{V_F}{\sigma_s} \bigr)^{1/s}.
\ea
$$
\EE

Let us also prove the following useful \textit{lifting} property.

\BL For any $q \geq 2$ and any $2 \leq s \leq q$, we have
\beq \label{Lifting}
\ba{rcl}
\bigl(  \frac{D_q}{D} \bigr)^q & \leq & 
\bigl(  \frac{D_s}{D} \bigr)^s.
\ea
\eeq
\EL
\proof
Indeed,
$$
\ba{rcl}
D_q & = &
\sup\limits_{\substack{x, y \in \mathcal{F}_0 \\ x \not= y}} 
\Bigl\{
\;
\frac{\|x - y\|}{\xi_F(x, y)^{1/q}}
\;
\Bigr\}
\;\; = \;\;
\sup\limits_{\substack{x, y \in \mathcal{F}_0 \\ x \not= y}} 
\Bigl\{
\;
\|x - y\|^{1 - \frac{s}{q}}
\bigr( \frac{\|x - y\|}{\xi_F(x, y)^{1/s}} \bigr)^{\frac{s}{q}}
\;
\Bigr\} \\
\\
& \leq & 
\sup\limits_{\substack{x, y \in \mathcal{F}_0 \\ x \not= y}} 
\Bigl\{
\;
\|x - y\|^{1 - \frac{s}{q}}
\;
\Bigr\} D_s^{\frac{s}{q}}
\;\; = \;\;
\bigr( D^{q - s} D_s^s \bigr)^{1/q}. \QR
\ea
$$

An immediate consequence of definition \eqref{RelSize}
comes from the Mean Value Theorem for convex functions:
$$
\ba{rcl}
F(y) & = & F(x) +\la G_x, y - x \ra + \int\limits_{0}^1
\frac{1}{\tau}
\la G_{x + \tau (y - x)} - G_x, \tau(y - x) \ra d\tau \\
\\
& \overset{\eqref{RelSize}}{\geq} & 
F(x) +\la G_x, y - x \ra +  
\frac{1}{s} V_F \Bigl( \frac{1}{D_s} \|y - x\|  \Bigr)^s,
\qquad \forall x, y \in \mathcal{F}_0.
\ea
$$ 
Hence, when $D_s < +\infty$, by minimizing the left and right hand 
sides with respect to $y$ independently, we get
\beq \label{UConvRes}
\ba{rcl}
F_* & \geq & 
F(x) - \frac{s - 1}{s}  \Bigl( \frac{D_s^s\| G_x \|_{*}^{s}  }{V_F} \Bigr)^{\frac{1}{s - 1}}
\;
\Leftrightarrow
\;\;
\frac{s - 1}{s}
\Bigl(  \frac{D_s \|G_x\|_{*}}{V_F} \Bigr)^{\frac{s}{s - 1}}
\; \geq \; 
\frac{F(x) - F_*}{V_F}.
\ea
\eeq

Let us introduce a formal assumption on our objective.

\BAS\label{AssumptionUConv} For some $s \geq 2$, we have $D_s < +\infty$.
\EAS
If $s = 2$, this assumption implies strong convexity.
If $s = \infty$, it means that the set $\mathcal{F}_0$ is bounded.

Since we define the relative size $D_s$ for the entire composite objective,
parameter $s \geq 2$ is not necessarily consistent 
with the degree of smoothness $q \in [2, 4]$.
Let us analyze the convergence rate of our method for different ranges of these parameters. 
We start with establishing the bound on the functional progress
during one iteration.

\BL \label{LemmaOneStepSConv}
For any $q \in [2, 4]$ and $s \in [2, \infty]$, we have
\beq \label{OneStepSConv}
\ba{rcl}
\frac{1}{(\gamma - 1) F_{k + 1}^{\gamma - 1}}
- \frac{1}{(\gamma - 1) F_k^{\gamma - 1}}
& \geq & \omega_{q, s} \bigl( \frac{g_{k + 1}}{g_k} \bigr)^2, 
\qquad k \geq 0,
\ea
\eeq
where $\gamma \Def \frac{q(s - 1)}{s(q - 1)} \in [\frac{2}{3}, 2]$
and $\omega_{q, s} \Def \frac{1}{16} \bigl(  \frac{s}{s - 1} \bigr)^{\gamma} 
\bigl(  \frac{V_F^{q/s}}{6 M_q D_s^{q}} \bigr)^{1 / (q - 1)}$.
\EL
\BR
Since $\lim\limits_{\alpha \to 0} \frac{x^{\alpha} - y^{\alpha}}{\alpha}
= \lim\limits_{\alpha \to 0} \frac{e^{\alpha \ln x} - e^{\alpha \ln y}}{\alpha} = \ln \frac{x}{y}$,
for $\gamma = 1 \Leftrightarrow s = q$, we treat the left-hand side of \eqref{OneStepSConv}
as its limit, which gives
$$
\ba{rcl}
\lim\limits_{\gamma \to 1}
\Bigl[ \frac{1}{(\gamma - 1) F_{k + 1}^{\gamma - 1}}
- \frac{1}{(\gamma - 1) F_k^{\gamma - 1}}  \Bigr] 
& = & \ln \frac{F_k}{F_{k + 1}}.
\ea
$$
\ER
\proof
For one step of the method, we have
\beq \label{OneStepSconv2}
\ba{rcl}
F_k - F_{k + 1} & \overset{\eqref{AdaptiveOneStep}}{\geq} &
\frac{1}{16 (6 M_q)^{1/(q - 1)}}
\bigl( \frac{g_{k + 1}}{g_k} \bigr)^{2} g_k^{\frac{q}{q - 1}} \\
\\
& \overset{\eqref{UConvRes}}{\geq} &
\frac{1}{16 (6 M_q)^{1/(q - 1)}}
\bigl(  \frac{g_{k + 1}}{g_k} \bigr)^2
\Bigl( \frac{V_F^{1/s}}{D_s}  \Bigr)^{\frac{q}{q - 1}}
\bigl(  \frac{s}{s - 1} F_k \bigr)^{\frac{q(s - 1)}{s(q - 1)}} \\
\\
& = & 
\omega_{q, s} \bigl(  \frac{g_{k + 1}}{g_k} \bigr)^2 F_k^{\gamma},
\ea
\eeq

First, let us consider the case $s \geq q$.
Then, $\gamma \in [1, 2]$. Using concavity of $y(x) = x^{\gamma - 1}, x \geq 0$,
and monotonicity of $\{ F_k \}_{k \geq 0}$,
we obtain
$$
\ba{rcl}
\frac{1}{ (\gamma - 1)F_{k + 1}^{\gamma - 1}}
- \frac{1}{(\gamma - 1) F_k^{\gamma - 1}}
& = &
\frac{
	F_k^{\gamma - 1} - F_{k + 1}^{\gamma - 1}
}{(\gamma - 1) F_{k + 1}^{\gamma - 1} F_{k}^{\gamma - 1}}
\;\; \geq \;\;
\frac{F_{k} - F_{k + 1}}{ F_{k + 1}^{\gamma - 1} F_k }
\;\; \overset{\eqref{OneStepSconv2}}{\geq} \;\;
\omega_{q, s} \bigl(  \frac{g_{k + 1}}{g_k}  \bigr)^2. 
\ea
$$

When $2 \leq s < q$, we have $\gamma < 1$. In this case,
we can use concavity of $y(x) = x^{1 - \gamma}, x \geq 0$.
This yields
$$
\ba{rcl}
\frac{1}{ (\gamma - 1)F_{k + 1}^{\gamma - 1}}
- \frac{1}{(\gamma - 1) F_k^{\gamma - 1}}
& = &
\frac{F_k^{1 - \gamma}}{1 - \gamma}
- \frac{F_{k + 1}^{1 - \gamma}}{1 - \gamma}
\;\; \geq \;\;
\frac{F_k - F_{k + 1}}{F_k^{\gamma}}
\;\; \overset{\eqref{OneStepSconv2}}{\geq} \;\;
\omega_{q, s} \bigl(  \frac{g_{k + 1}}{g_k}  \bigr)^2. \quad \QR
\ea
$$

Inequality~\eqref{OneStepSConv} provides us 
with a \textit{continuous} in $\gamma$ characterization
of the global behavior of the method.
We are ready to describe its global complexity. 
Let us start with the case $s \geq q$.

\BT \label{TheoremGlobalUConvex}
Let $M_q < +\infty$ for some $q \in [2, 4]$,
and  $D_s < +\infty$ for $\boxed{s \geq q}$.
Assume that the initial value $H_0$ satisfies~\eqref{H_0_Bound}
and that the functional residual for all iterations 
$\{ x_i \}_{i = 0}^{k}$ of Algorithm~\ref{alg:AdaptiveNewton}
is big enough:
\beq \label{ResBoundedBelowUConv}
\ba{rcl}
F_i & \Def & F(x_i) - F_* \;\; \geq \;\; \varepsilon,
\ea
\eeq
with some $\varepsilon > 0$. Then,
\beq \label{UConvComplexity}
\ba{rcl}
k & \leq & 
16 \bigl( \frac{s - 1}{s} \bigr)^{\frac{q(s - 1)}{s(q - 1)}}
\Bigl(  \frac{6 M_q D_s^{q} }{V_F^{q/s}} \Bigr)^{\frac{1}{q - 1}}
\frac{s(q - 1)}{s - q}
\Bigl[
\varepsilon^{-\frac{s - q}{s(q - 1)}} - F_0^{-\frac{s - q}{s(q - 1)}}
\Bigr]
 \!+\! 2 \ln \frac{g_0 D}{\varepsilon}.
\ea
\eeq
For $\boxed{s = q}$, we treat the right-hand side of \eqref{UConvComplexity}
as its limit, which gives
the linear convergence rate:
\beq \label{UconvLinearComplexity}
\ba{rcl}
k & \overset{\eqref{UConvComplexity}}{\leq} &
16 \bigl( \frac{q - 1}{q}  \bigr)
\Bigl( \frac{6M_q D_q^q}{V_F}  \Bigr)^{\frac{1}{q - 1}}
\ln \frac{F_0}{\varepsilon}
+ 2 \ln \frac{g_0 D}{\varepsilon}.
\ea
\eeq
\ET
\proof
Telescoping bound~\eqref{OneStepSConv} and using the inequality for 
arithmetic and geometric means, we get
$$
\ba{rcl}
\frac{1}{(\gamma - 1) F_k^{\gamma - 1}}
- \frac{1}{(\gamma - 1) F_0^{\gamma - 1}}
& \overset{\eqref{OneStepSConv}}{\geq} & 
\omega_{q,s} \sum\limits_{i = 0}^{k - 1} 
\bigl( \frac{g_{i + 1}}{g_i} \bigr)^2
\;\; \geq \;\;
k \omega_{q, s} \bigl( \frac{g_k}{g_0} \bigr)^{2/k} \\
\\
& \overset{\eqref{ConvexityFD}}{\geq} &
k \omega_{q, s} \bigl( \frac{F_k}{g_0 D} \bigr)^{2/k}
\;\; = \;\;
k \omega_{q, s} \cdot \exp\Bigl( -\frac{2}{k} \ln \frac{g_0 D}{F_k}  \Bigr) \\
\\ 
& \geq &
k \omega_{q, s} \cdot \Bigl(1 - \frac{2}{k} \ln \frac{g_0 D}{F_k} 	\Bigr) 
\;\; = \;\;
k \omega_{q, s} - 2 \omega_{q, s} \ln\frac{g_0D}{F_k}.
\ea
$$
Therefore, 
$$
\ba{rcl}
k & \leq & \frac{1}{ \omega_{q, s}(\gamma - 1)}
	\Bigl[ 
	\frac{1}{F_k^{\gamma - 1}} - \frac{1}{F_0^{\gamma - 1}}
	\Bigr]
+ 2\ln\frac{g_0 D}{F_k}.
\ea
$$
Substituting the bound \eqref{ResBoundedBelowUConv} for $F_k$,
and using the definitions of $\gamma$ and $\omega_{q, s}$,
we obtain \eqref{UConvComplexity}.
\qed

Note that for $s = q = 2$ (strongly convex functions
with bounded variation of the Hessian),
inequality \eqref{UconvLinearComplexity} implies
global linear rate. Thus, it covers the standard problem class for the Gradient Methods.

Let us consider now $2 \leq s < q$.
According to \eqref{Lifting}, we have
$D_q^q \leq D^{q - s} D_s^s$.
Substituting this bound into~\eqref{UconvLinearComplexity},
we would get the complexity estimate for this case.
However, it would only give us a \textit{linear} rate.
As we will see, if $2 \leq s < q$, the method
has a \textit{superlinear} convergence.

\BT \label{TheoremGlobalSuperlinear}
Let $M_q < +\infty$ for some $q \in [2, 4]$,
and $D_s < +\infty$ for $\boxed{2 \leq s < q}$.
Assume that the initial value $H_0$ satisfies \eqref{H_0_Bound},
and for all iterations $\{ x_i \}_{i = 0}^k$
of Algorithm~\ref{alg:AdaptiveNewton}, we have
\beq \label{GlobalSuperGradBig}
\ba{rcl}
F_i & \Def & F(x_i) - F_* \;\; \geq \;\; \varepsilon,
\qquad
g_i \;\; \Def \;\; \| F'(x_i) \|_{*} \;\; \geq \;\; \delta,
\ea
\eeq
for some $\varepsilon > 0$ and $\delta > 0$. Then,
\beq \label{GlobalSuperlinear}
\ba{rcl}
k & \leq & 
\frac{16 q}{s} \bigl( \frac{s - 1}{s} \bigr)^{\frac{s - 1}{q - 1}}
\Bigl(   6 M_q  \frac{D_s^s D^{q - s}}{V_F}   \Bigr)^{\frac{1}{q - 1}} \\[10pt]
& \; & \; 
\times \;\;
\frac{s(q - 1)}{q - s} \Bigl[ 
1 -   \frac{s}{q} \Bigl( \bigl( \frac{s - 1}{s} \bigr)^{s - 1} 
\frac{ D_s^s}{D^s V_F} \, \varepsilon \Bigr)^{\frac{q - s}{s(q - 1)}}
\Bigr]
\; + \; 2 \ln\frac{g_0}{\delta}.
\ea
\eeq
For $\boxed{s = q}$, we treat the right-hand side of \eqref{GlobalSuperlinear}
as its limit, which gives the linear rate of convergence:
$$
\ba{rcl}
k & \overset{\eqref{GlobalSuperlinear}}{\leq} & 
16 \bigl(  \frac{q - 1}{q} \bigr)
\Bigl( \frac{6 M_q D_q^q}{V_F}  \Bigr)^{\frac{1}{q - 1}}
\ln\Bigl[  \bigl(  \frac{q}{q - 1} \bigr)^{q - 1}  \frac{D^q V_F}{D_q^q \varepsilon} \Bigr]
\; + \;
2 \ln\frac{g_0}{\delta}.
\ea
$$
\ET
\proof
We split iterations of the method into two consecutive stages
$k = m + n$. During the first $m$ iterations, we use bound~\eqref{ConvexTelescopedStep}, that provides us
with the guarantee
\beq \label{GSupPart1}
\ba{cl}
& C_q \bigl(\frac{ 1 }{F_m} \bigr)^{\frac{1}{q - 1}}
\;\; \geq \;\;
C_q \Bigl[   \bigl(\frac{1}{F_m}\bigr)^{\frac{1}{q - 1}} - 
      \bigl(\frac{1}{F_0}\bigr)^{\frac{1}{q - 1}}  \Bigr]
\;\; \overset{\eqref{ConvexTelescopedStep}}{\geq} \;\;
m \bigl( \frac{g_m}{g_0} \bigr)^{\frac{2}{m}} \\
\\
& = \;\;
m \exp\bigl( -\frac{2}{m} \ln \frac{g_0}{g_m} \bigr)
\;\; \geq \;\;
m - 2 \ln \frac{g_0}{g_m}
\;\; = \;\;
m + 2 \ln \frac{g_m}{\delta} - 2\ln\frac{g_0}{\delta}.
\ea
\eeq
where $C_q \Def 16 (q - 1) (6 M_q)^{\frac{1}{q - 1}} D^{\frac{q}{q - 1}}$.
During the second stage, which lasts $n$ iterations, we have
$$
\ba{rcl}
\frac{1}{\omega_{q, s} (1 - \gamma)}
\bigl[ F_i^{1 - \gamma} - F_{i + 1}^{1 - \gamma}  \bigr]
& \overset{\eqref{OneStepSConv}}{\geq} &
\bigl(  \frac{g_{i + 1}}{g_i} \bigr)^{2},
\quad i = m, \ldots, k - 1.
\ea
$$
Note that $1 - \gamma = \frac{q - s}{s(q - 1)} > 0$.
Telescoping these inequalities, we get
\beq \label{GSupPart2}
\ba{cl}
& \frac{1}{\omega_{q, s}(1 - \gamma)}
\bigl[ F_m^{1 - \gamma} - \varepsilon^{1 - \gamma}  \bigr]
\;\; \overset{\eqref{GlobalSuperGradBig}}{\geq} \;\;
\frac{1}{\omega_{q, s}(1 - \gamma)}
\bigl[ F_m^{1 - \gamma} - F_k^{1 - \gamma}  \bigr] \\
\\
& \geq \;\;
\sum\limits_{i = m}^{k - 1} \bigl(  \frac{g_{i + 1}}{g_i} \bigr)^{2}
\;\; \geq \;\;
n \bigl( \frac{g_{k}}{g_{m}} \bigr)^{\frac{2}{n}}
\;\; \overset{\eqref{GlobalSuperGradBig}}{\geq} \;\;
n \bigl(  \frac{\delta}{g_m} \bigr)^{\frac{2}{n}}
\;\; \geq \;\;
n - 2 \ln \frac{g_m}{\delta}.
\ea
\eeq
Hence,
$$
\ba{rcl}
k & = & m + n \;\; \overset{\eqref{GSupPart1}, \eqref{GSupPart2}}{\leq} \;\;
C_q \bigl( \frac{1}{F_m} \bigr)^{\frac{1}{q - 1}}
+ \frac{1}{\omega_{q, s}(1 - \gamma)}
\bigl[ F_m^{1 - \gamma} - \varepsilon^{1 - \gamma}  \bigr] 
+ 2\ln\frac{g_0}{\delta}.
\ea
$$
The maximum of the right-hand side as a function of $F_m$ is attained at
$$
\ba{rcl}
\tau^{*} & = & 
\argmax\limits_{\tau > 0}
\Bigl\{ 
C_q \bigl(  \frac{1}{\tau} \bigr)^{\frac{1}{q - 1}}
+ \frac{1}{w_{q, s} (1 - \gamma)} \tau^{1 - \gamma} 
\Bigr\}
\;\; = \;\;
\Bigl( \frac{C_q \omega_{q,s}}{q - 1} \Bigr)^{\frac{s(q - 1)}{q}}.
\ea
$$
Substituting this value, we get
$$
\ba{rcl}
k & \leq & 
C_q^{\frac{q - s}{q}}
\bigl( \frac{q - 1}{\omega_{q,s}}  \bigr)^{\frac{s}{q}}
+ \frac{s(q - 1)}{q - s}
\Bigl[
\bigl(  \frac{C_q}{q - 1} \bigr)^{\frac{q - s}{s}}
\bigl( \frac{1}{\omega_{q,s}}  \bigr)^{\frac{s}{q}}
- \frac{1}{\omega_{q, s}} \varepsilon^{\frac{q - s}{s(q - 1)}}
\Bigr] 
+ 2\ln\frac{g_0}{\delta}.
\ea
$$
Using the definitions of $C_q$ and $\omega_{q, s}$ completes the proof. \qed

It remains to analyze local convergence of the method.
Denoting $r_k = \|x_k - x_{k + 1}\|$, we have
\beq \label{LocalBounds}
\ba{cl}
& g_k r_k \;\; \geq \;\; \la F'(x_k), x_k - x_{k + 1} \ra \\
\\
& \overset{\eqref{RelSize}}{\geq} \;
\la F'(x_{k + 1}), x_k - x_{k + 1} \ra + V_F \bigl(\frac{r_k}{D_s}\bigr)^s 
\;\; \geq \;\;
V_F \bigl(  \frac{r_k}{D_s} \bigr)^s.
\ea
\eeq
Hence,
\beq \label{Superlinear}
\ba{rcl}
 g_{k + 1} 
& \overset{\eqref{NewGradBound}}{\leq} &
4 \lambda_k r_k
\;\; \overset{\eqref{LocalBounds}}{\leq} \;\;
4 \lambda_k \bigl( \frac{g_k D_s^s}{ V_F } \bigr)^{\frac{1}{s - 1}}  \\
\\
& \overset{\eqref{Lambda_k_Bound}}{\leq} &
16 \bigl( 6 M_q \bigr)^{\frac{1}{q - 1}} 
\bigl(  \frac{D_s^s}{V_F} \bigr)^{\frac{1}{s - 1}}
g_k^{ \frac{1}{s - 1} + \frac{q - 2}{q - 1} } 
\;\; \equiv \;\;
\bigl(\frac{1}{\Delta}\bigr)^{\zeta}  g_k^{1 + \zeta},
\ea 
\eeq
where $\zeta \Def \frac{q - s}{(s - 1)(q - 1)}$
and 
$
\Delta
 \Def  \Bigl[
\frac{V_F^{ 1 / (s - 1) } }{ 16 D_s^{s/(s - 1)} (6 M_q)^{1 / (q - 1)} }
\Bigr]^{\! \frac{1}{\zeta}}.
$
So, we have just proved the local superlinear convergence of power $1 + \zeta$.
\BT \label{TheoremLocalSuperlinear}
Let $\boxed{2 \leq s < q}$.
Assume the gradient is small enough:
\beq \label{SuperlinearRegion}
\ba{rcl}
g_0 & \in & \mathcal{Q}
\;\; \Def \;\;
\Bigl\{
g \; : \; 
\| g \|_{*} \; \leq \; \frac{1}{e}\Delta
 \Bigr\}.
\ea
\eeq
Then, for any $\delta > 0$, by doing
\beq \label{SuperlinearRate}
\ba{rcl}
k & = & 
\Bigl\lceil \frac{1}{\ln(1 + \zeta)} \ln\ln \frac{\Delta}{\delta} \Bigr\rceil
\ea
\eeq
iterations, we have $g_k \leq \delta$.
\ET
\proof
Dividing both sides of~\eqref{Superlinear} by $\Delta$,
we get: $\frac{g_{k + 1}}{\Delta} \leq \bigl( \frac{g_k}{\Delta} \bigr)^{1 + \zeta}$. Hence,
$$
\ba{rcl}
\frac{g_k}{\Delta} & \leq & 
\bigl( \frac{g_{k - 1}}{\Delta} \bigr)^{1 + \zeta} 
\;\; \leq \;\; \ldots \;\; \leq \;\;
\bigl(  \frac{g_0}{\Delta} \bigr)^{(1 + \zeta)^k}
\;\; \overset{\eqref{SuperlinearRegion}}{\leq} \;\;
\bigl( \frac{1}{e} \bigr)^{(1 + \zeta)^k}.
\ea
$$
After \eqref{SuperlinearRate} iterations, we 
ensure $\bigl( \frac{1}{e} \bigr)^{(1 + \zeta)^k} \leq \frac{\delta}{\Delta}$,
which completes the proof.
\qed

For example, for $q = 3$ (Lipschitz Hessian) and $s = 2$ (strongly convex functions), 
inequality~\eqref{Superlinear} implies
$$
\ba{rcl}
g_{k + 1} & \leq & O\Bigl(  \frac{M_3^{1/2} D_2^2 }{V_F}  g_k^{3/2}  \Bigr)
\;\; = \;\; 
O\Bigl(  \frac{M_3^{1/2}}{\mu} g_k^{3/2} \Bigr),
\ea
$$	 
where $\mu > 0$ is the parameter of strong convexity.
Note that this is slightly worse than the local \textit{quadratic} convergence
of the pure Newton Method \cite{nesterov2018lectures}.
However, it seems to be a reasonable price for the universality.

For $q = 4$ (Lipschitz third derivative) and $s = 2$, we get
$$
\ba{rcl}
g_{k + 1} & \leq & O\Bigl(  \frac{M_4^{1/3} D_2^2}{V_F}  g_k^{5/3}  \Bigr)
\;\; = \;\;
O\Bigl(  \frac{M_4^{1/3}}{\mu} g_k^{5/3} \Bigr).
\ea
$$	 

Note that now the target accuracy $\delta$ enters into~\eqref{SuperlinearRate} under \textit{two logarithms}.
This is a very fast convergence, and for all practical applications it is enough to do 
only a constant number of steps, after reaching~\eqref{SuperlinearRegion}.

\section{Numerical Experiments}
\label{SectionExperiments}

\paragraph{Polytope Feasibility.}
We model the
problem of finding a feasible point $x^{*} \in \mathcal{P}$ of a polytope 
$\mathcal{P} = \bigl\{  x \in \R^n : \; \la a_i, x \ra \leq b_i, \; 1 \leq i \leq m  \bigr\}$,
by the following minimization objective:
$$
\ba{rcl}
\min\limits_{x \in \R^n} \Bigl[ \;f(x) 
& := &
\sum\limits_{i = 1}^m \bigl(  \la a_i, x \ra - b_i \bigr)_+^p \;\Bigr],
\ea
$$
where $(t)_{+} \Def \max\{ 0, t \}$ is positive slicing and $p \geq 2$ is our parameter.

We compare our method for $\alpha = \frac{2}{3}$ and $\alpha = 1$ with 
the following algorithms:
Cubic Newton, Gradient Method, and Fast Gradient Method \cite{nesterov2018lectures}.
In all the methods we use adaptive estimation of the regularization parameter.

We generate data from random uniform distribution on $[-1, 1]$, and start the methods from $x_0 = (1, \ldots, 1)^\top \in \R^n$.
In many cases, we had degenerate Hessian at the initial point $\nabla^2 f(x_0)$, and so it is impossible 
to use a Damped Newton Method. The results are shown 
below\footnote{The source code can be found at
\url{https://github.com/doikov/super-newton/}}.

\begin{figure}[H]
	\centering
	\includegraphics[width=0.49\textwidth ]{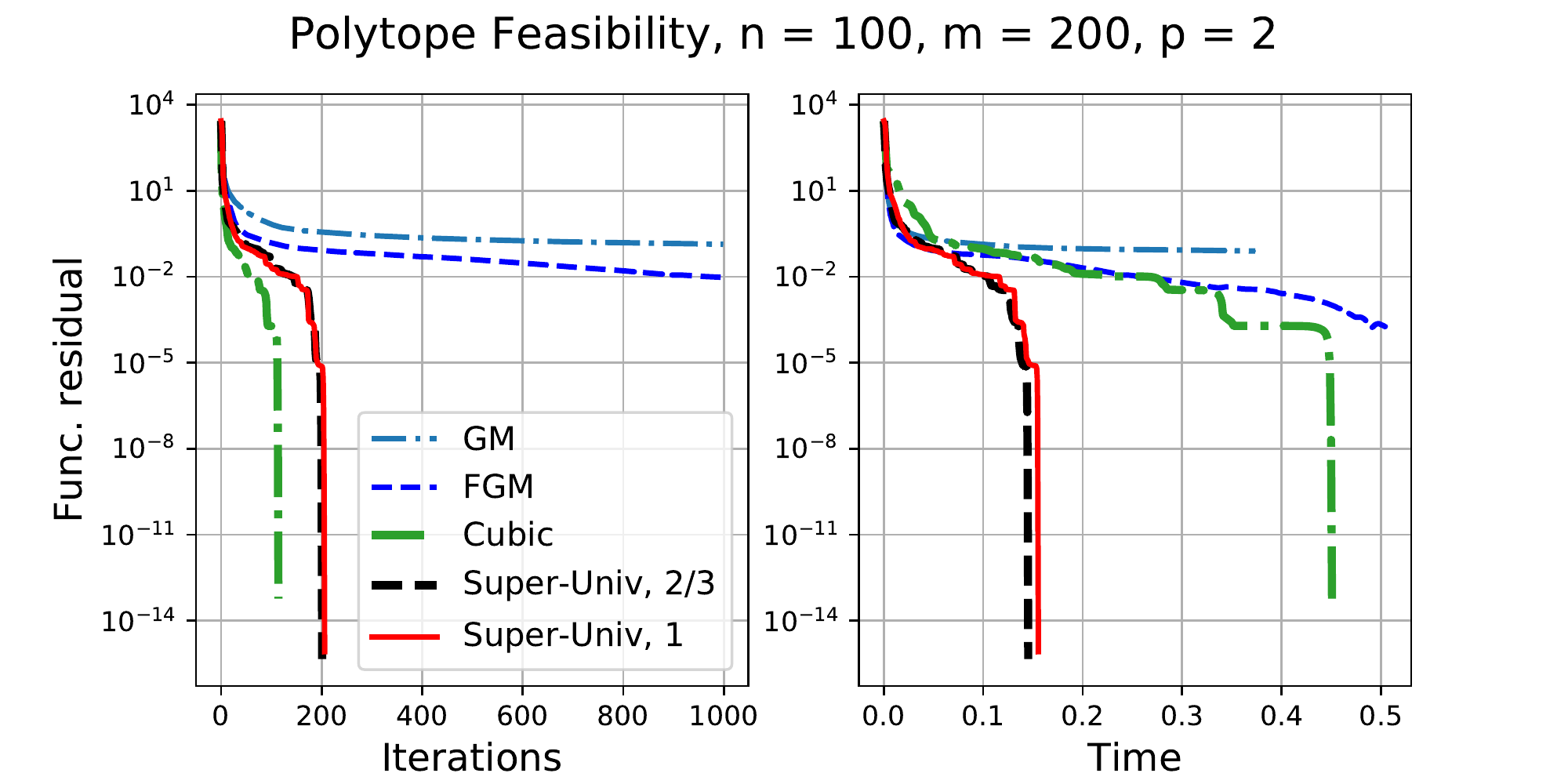}
	\includegraphics[width=0.49\textwidth ]{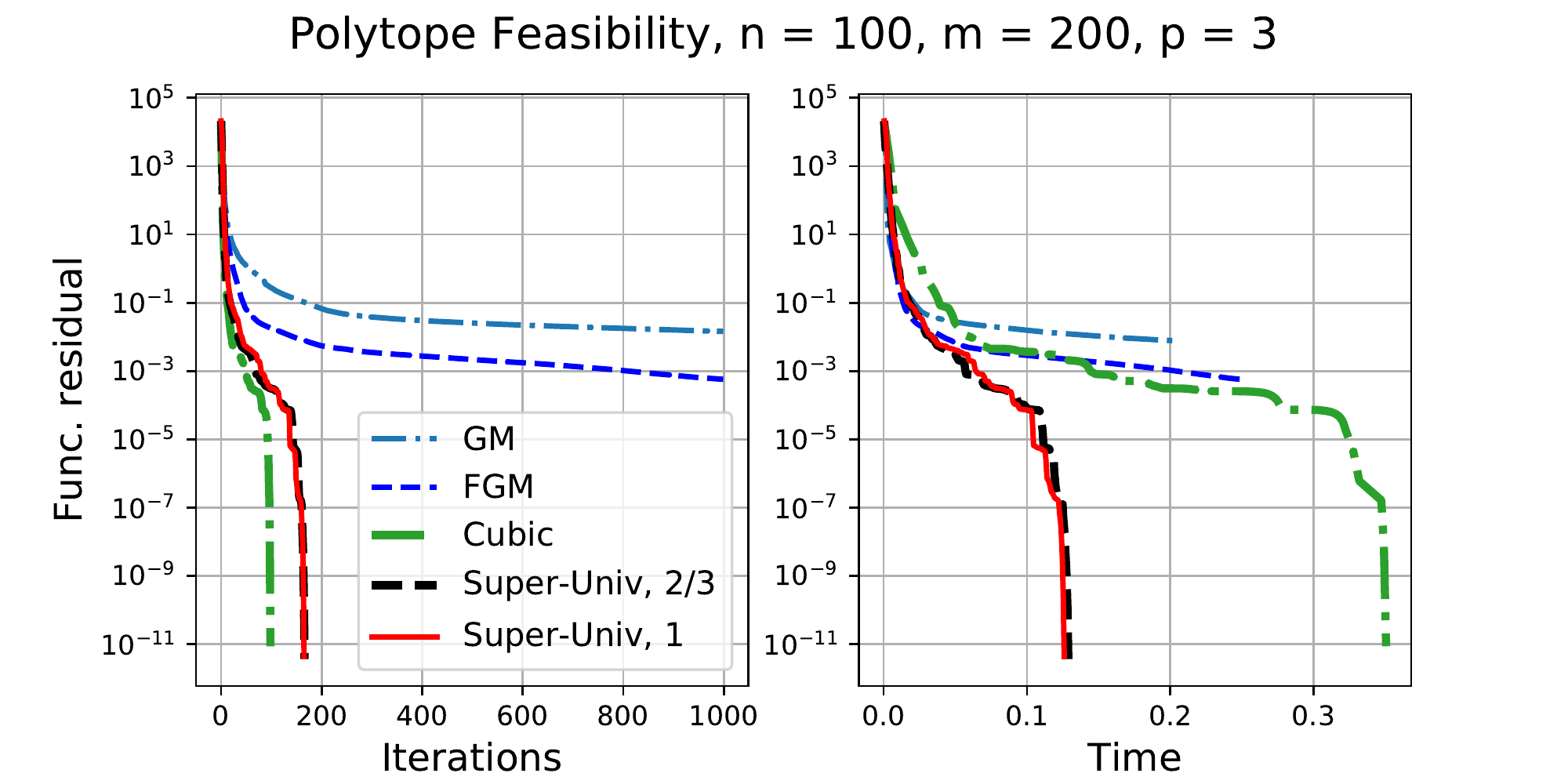}
	
	\includegraphics[width=0.49\textwidth ]{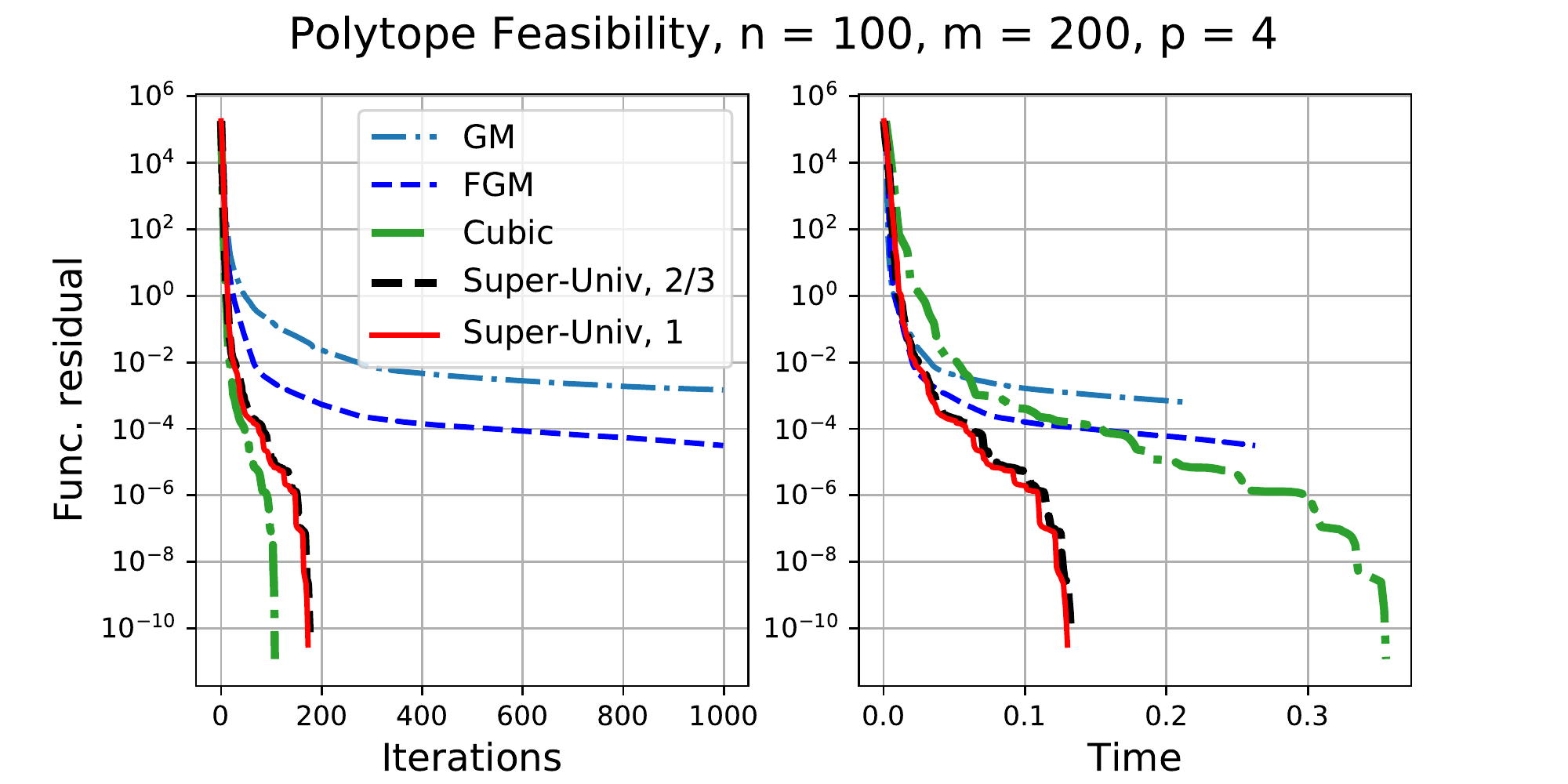}
	\includegraphics[width=0.49\textwidth ]{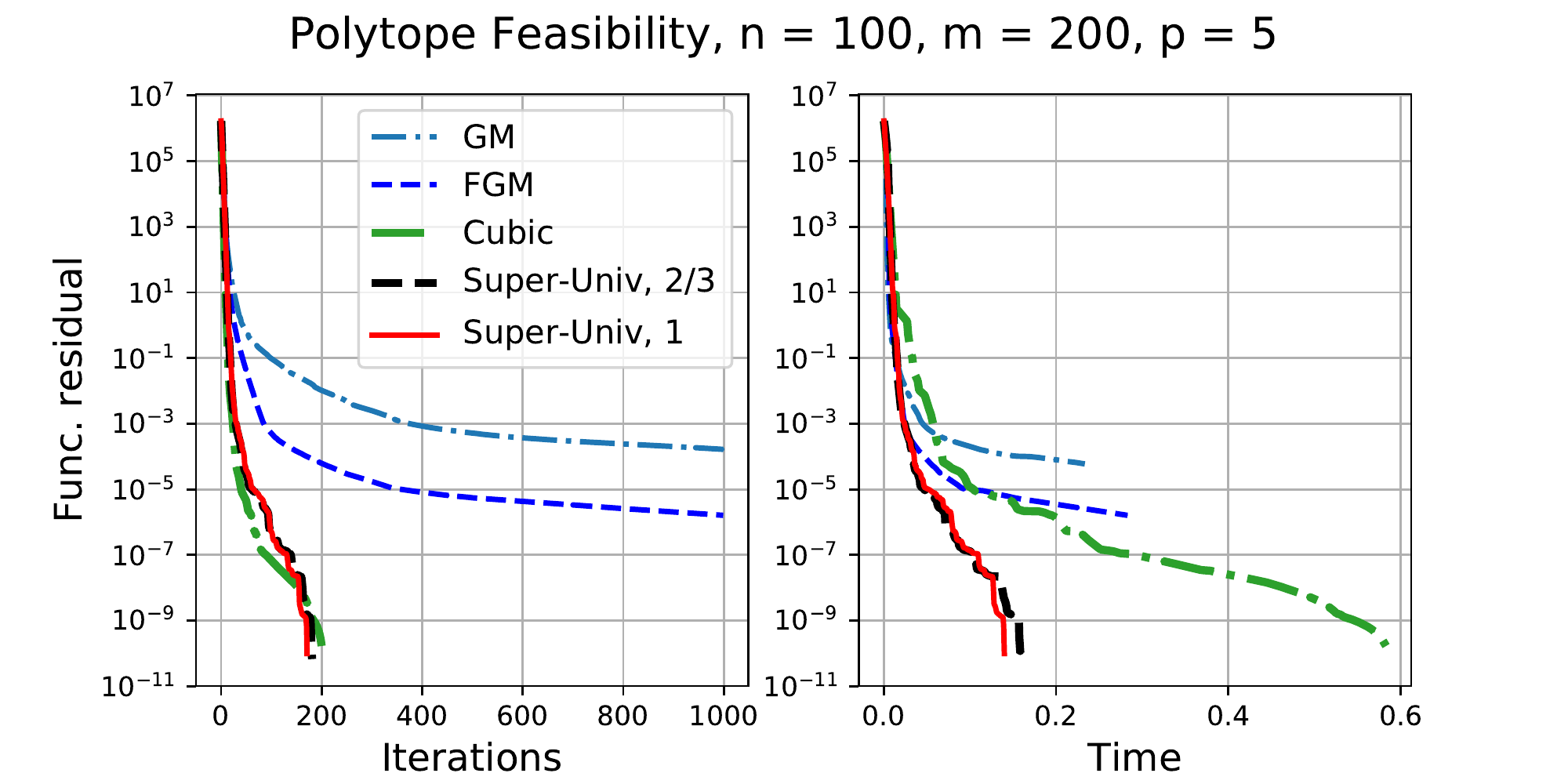}
\end{figure}

Thus, the second-order methods demonstrate extremely good performance
in terms of the number of iterations (oracle calls). 
The practical convergence of the Cubic Newton seems to be slightly better
than those with quadratic regularization. However, each iteration
of the latter methods is cheaper, which results in much better total computational time (see the second graph in each pair).

\paragraph{Soft Maximum.}
For $\mu > 0$,
consider the unconstrained minimization problem $\min_{x \in \R^n} f(x)$
with the following objective, 
$$
\ba{rcl}
f(x)  & := & \mu \ln \Bigl(  \; \sum\limits_{i = 1}^m 
\exp\bigl(  \frac{\la a_i, x \ra - b_i}{\mu}\bigr)
\, \Bigr)
\quad \approx \quad \max\limits_{1 \leq i \leq m} \bigl[ \la a_i, x \ra - b_i \bigr].
\ea
$$
The entries of vectors $a_1, \ldots, a_m \in\ \R^n$ and $b \in \R^m$
are generated randomly and independently from the uniform distribution on
$[-1, 1]$, and $\mu$ is a smoothing parameter.
We use the primal norm, with the following matrix:
$B = \sum_{i = 1}^m a_i a_i^\top$.
Then we have (Example 1.3.5 in \cite{doikov2021new}):
$$
\ba{rcl}
M_2 & \leq & \frac{1}{\mu}, \quad
M_4 \; \leq \; \frac{4}{\mu^3}
\quad \Rightarrow \quad
M_q \; \overset{\eqref{MqBound}}{\leq} \;
3 \cdot \frac{4^{\frac{q - 2}{2}}}{\mu^{q - 1}}
\; \leq \; \frac{12}{\mu^{q - 1}}, \quad \forall q \in [2, 4].
\ea
$$
Thus, for $\varepsilon > 0$, our method needs to do the following number of iterations:
$$
\ba{rcl}
k & \overset{\eqref{ConvexTotalComplexity}}{=} &
O\Bigl( \; \frac{1}{\mu} \inf\limits_{2\leq q \leq 4}  
\bigl(  \frac{D^q}{\varepsilon} \bigr)^{\frac{1}{q- 1}}  + \ln\frac{g_0 D}{\varepsilon} \; \Bigr)
\ea
$$

By an appropriate shifting of all vectors $\{ a_i \}_{i = 1}^m$, 
we can ensure $\nabla f(0) = 0$, placing the optimum to the origin.
We run Algorithm~\ref{alg:AdaptiveNewton}
with different values of the gradient power $\alpha = 0, \frac{1}{2}, \frac{2}{3}, 1$.
The results are presented below.

\begin{figure}[H]
	\centering
	\includegraphics[width=0.48\textwidth ]{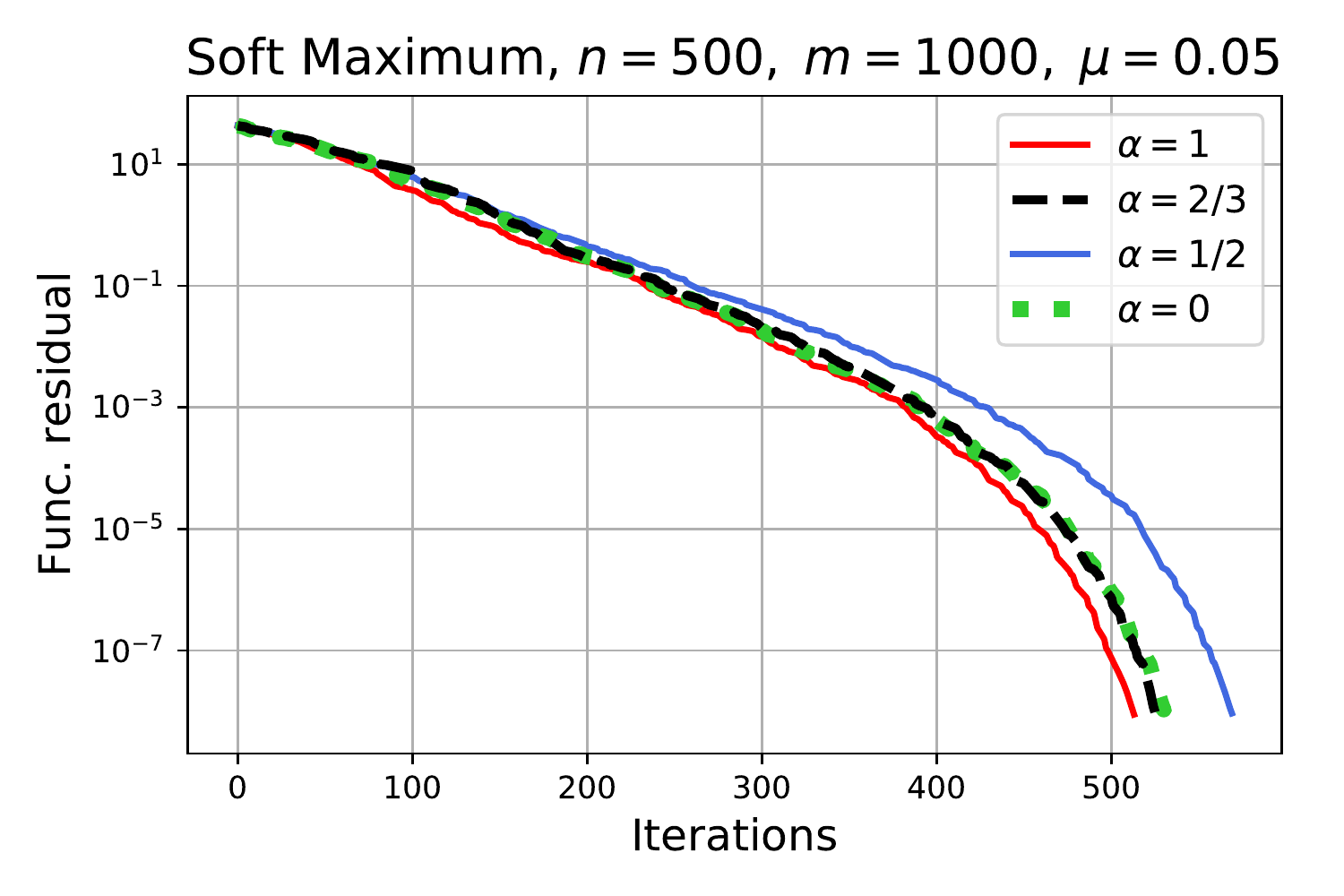}
	\includegraphics[width=0.48\textwidth ]{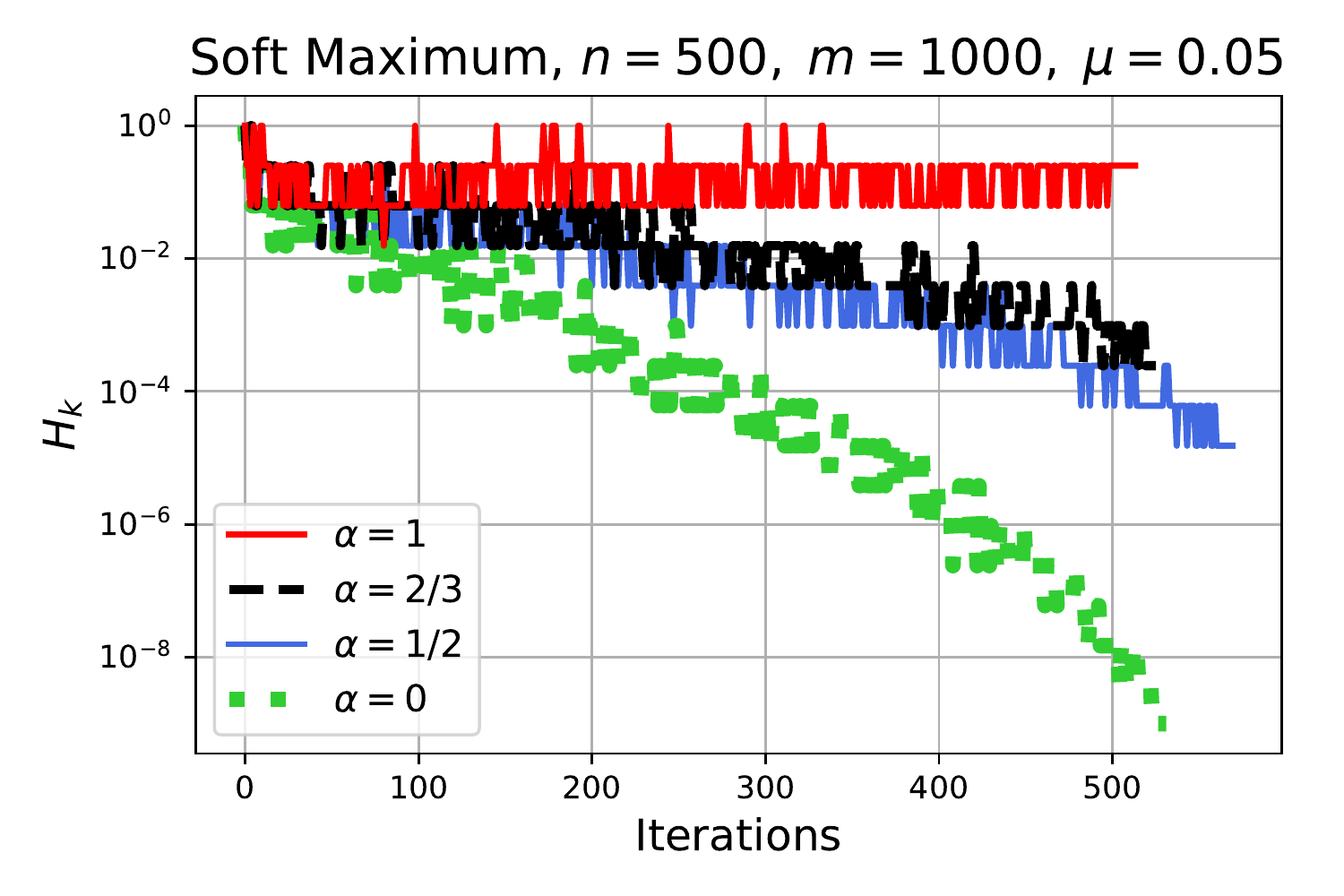}
\end{figure}

The method shows a robust behaviour 
in terms of the dependence on $\alpha$ (left graph). However, the choice $\alpha := 1$
enforces a more stable range for the regularization parameter $H_k$
adjusted by the adaptive search (right graph).

\paragraph{Worst Instances.}
In this experiment, we apply our methods to 
unconstrained minimization of the following objective,
$$
\ba{rcl}
f(x) & := &  \frac{1}{q} \sum\limits_{i = 1}^{n - 1}| x^{(i)} - x^{(i + 1)}|^q
+ \frac{1}{q}|x^{(n)}|^q, \qquad x \in \R^n,
\ea
$$
where $q \geq 2$ is a parameter. Note that the structure of this objective
is very similar to the worst-case function from lower bounds
for high-order methods \cite{nesterov2019implementable},
and there is a bound for the smoothness constants: $M_q \leq 2^q (q!)$.
We compare the method with fixed constants of regularization and
super-universal methods. The results are shown below.

\begin{figure}[H]
	\centering
	\includegraphics[width=0.40\textwidth ]{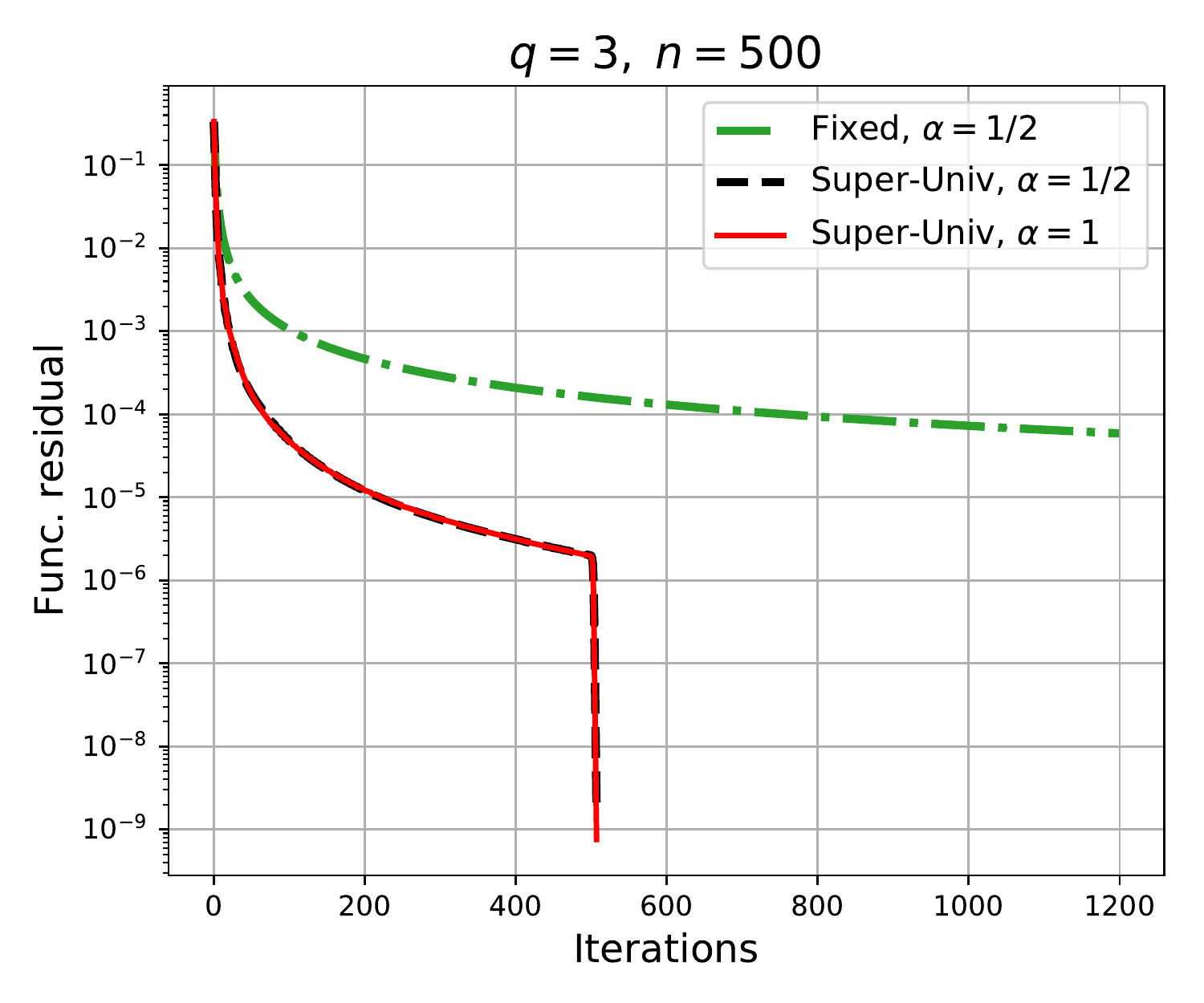}
	\includegraphics[width=0.40\textwidth ]{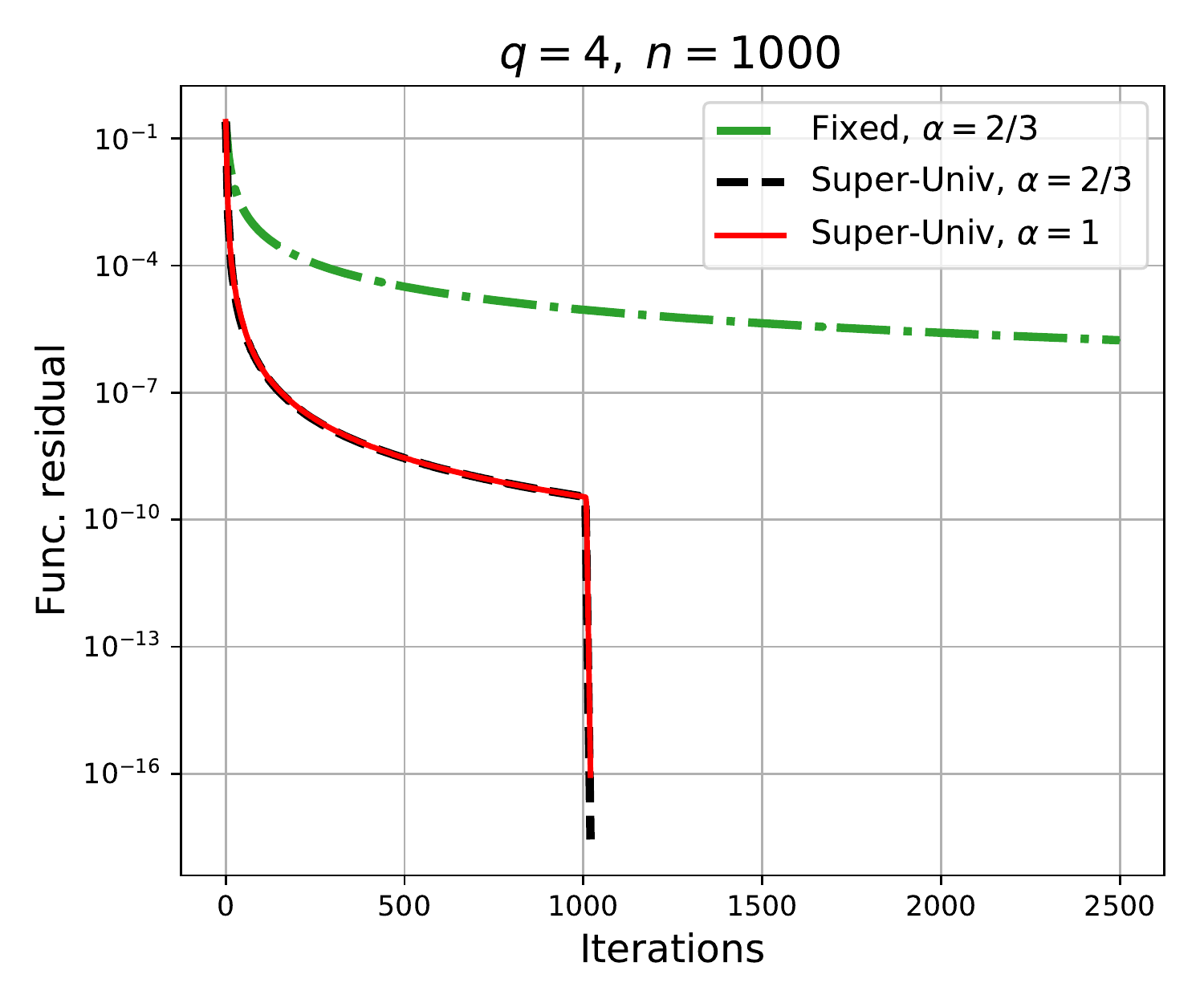}
\end{figure}
 
 We observe a switching point in the behaviour when the number of iterations
 reaches the dimensionality of the problem.
 In the case of super-universal methods, 
 the rate becomes superlinear after this moment,
 and any desirable accuracy can be achieved just in few extra steps.

\section{Discussion}
\label{SectionDiscussion}

In this paper, we have developed and analyzed the Super-Universal Newton Method
based on regularization of the second-order model by the square of Euclidean norm.
The regularization parameter is proportional to a power of the gradient norm.
Each step of our method is easily computable, employing
in the unconstrained case just the standard matrix inversion.

We have proved that using a simple adaptive search procedure in each iteration,
the method has a universal global convergence rate among problem classes
with H\"older continuous second or third derivatives.
If the problem is uniformly convex, the method automatically switches between
sublinear, linear, and superlinear rates, adjusting to the best possible
problem class.

A natural extension of our results would be development of \textit{accelerated} 
super-universal schemes (see \cite{nesterov2008accelerating,monteiro2013accelerated,grapiglia2019accelerated,grapiglia2020tensor,gasnikov2019near,doikov2020contracting,nesterov2021inexact,kovalev2022first,carmon2022optimal}
for the line of works on accelerated second- and high-order methods
matching the corresponding lower bounds \cite{arjevani2019oracle,nesterov2018lectures,grapiglia2020tensor}).
One of the major obstacles remains to be the sensitivity of accelerated methods
to the parameters of a problem class.
In addition, these methods usually require knowledge of the constant of strong/uniform convexity.
For practical applications, it is also crucial for a second-order method
to have a \textit{superlinear} convergence (at least locally),
which is missing for most of the accelerated schemes.

Another important direction is the creation of methods that are suitable for non-Euclidean geometry.
In our method we fix the Euclidean norm as a regularizer,
while it is also possible to use for that a contraction of the feasible domain, leading to 
affine-invariant \textit{contracting-point} methods \cite{doikov2022affine},
or an appropriate Bregman divergence \cite{doikov2021gradient}
(see also \cite{bauschke2016descent,lu2018relatively} for the framework of \textit{relative smoothness}).

For solving large-scale problems, our method can be equipped with
modern stochastic techniques \cite{rodomanov2016superlinearly,doikov2018randomized,kovalev2019stochastic,hanzely2020stochastic}
which are able to keep versatile convergence guarantees. Another potential way to make the methods more applicable to high-dimensional objectives is to consider quasi-Newton updates, which at the moment seems to be very challenging due to the lack of theoretical results on their global behavior.

\bibliographystyle{plain}
\bibliography{bibliography}

\begin{thebibliography}{10}

\bibitem{arjevani2019oracle}
Yossi Arjevani, Ohad Shamir, and Ron Shiff.
\newblock Oracle complexity of second-order methods for smooth convex
  optimization.
\newblock {\em Mathematical Programming}, 178(1-2):327--360, 2019.

\bibitem{baes2009estimate}
Michel Baes.
\newblock Estimate sequence methods: extensions and approximations.
\newblock {\em Institute for Operations Research, ETH, Z{\"u}rich,
  Switzerland}, 2009.

\bibitem{bauschke2016descent}
Heinz~H. Bauschke, J{\'e}r{\^o}me Bolte, and Marc Teboulle.
\newblock A descent lemma beyond {L}ipschitz gradient continuity: first-order
  methods revisited and applications.
\newblock {\em Mathematics of Operations Research}, 42(2):330--348, 2016.

\bibitem{birgin2017worst}
Ernesto~G. Birgin, J.~L. Gardenghi, Jos{\'e}~Mario Mart{\'\i}nez,
  Sandra~Augusta Santos, and Philippe~L. Toint.
\newblock Worst-case evaluation complexity for unconstrained nonlinear
  optimization using high-order regularized models.
\newblock {\em Mathematical Programming}, 163(1-2):359--368, 2017.

\bibitem{carmon2022optimal}
Yair Carmon, Danielle Hausler, Arun Jambulapati, Yujia Jin, and Aaron Sidford.
\newblock Optimal and adaptive {M}onteiro-{S}vaiter acceleration.
\newblock {\em arXiv preprint arXiv:2205.15371}, 2022.

\bibitem{cartis2011adaptive1}
Coralia Cartis, Nicholas I.~M. Gould, and Philippe~L. Toint.
\newblock Adaptive cubic regularisation methods for unconstrained optimization.
  {P}art {I}: motivation, convergence and numerical results.
\newblock {\em Mathematical Programming}, 127(2):245--295, 2011.

\bibitem{cartis2011adaptive2}
Coralia Cartis, Nicholas I.~M. Gould, and Philippe~L. Toint.
\newblock Adaptive cubic regularisation methods for unconstrained optimization.
  {P}art {II}: worst-case function-and derivative-evaluation complexity.
\newblock {\em Mathematical programming}, 130(2):295--319, 2011.

\bibitem{cartis2020sharp}
Coralia Cartis, Nicholas I.~M. Gould, and Philippe~L. Toint.
\newblock Sharp worst-case evaluation complexity bounds for arbitrary-order
  nonconvex optimization with inexpensive constraints.
\newblock {\em SIAM Journal on Optimization}, 30(1):513--541, 2020.

\bibitem{conn2000trust}
Andrew~R. Conn, Nicholas I.~M. Gould, and Philippe~L. Toint.
\newblock {\em Trust region methods}.
\newblock SIAM, 2000.

\bibitem{doikov2021new}
Nikita Doikov.
\newblock {\em New second-order and tensor methods in Convex Optimization}.
\newblock PhD thesis, Universit{\'e} catholique de Louvain, 2021.

\bibitem{doikov2020contracting}
Nikita Doikov and Yurii Nesterov.
\newblock Contracting proximal methods for smooth convex optimization.
\newblock {\em SIAM Journal on Optimization}, 30(4):3146--3169, 2020.

\bibitem{doikov2021gradient}
Nikita Doikov and Yurii Nesterov.
\newblock Gradient regularization of {N}ewton method with {B}regman distances.
\newblock {\em arXiv preprint arXiv:2112.02952}, 2021.

\bibitem{doikov2021local}
Nikita Doikov and Yurii Nesterov.
\newblock Local convergence of tensor methods.
\newblock {\em Mathematical Programming}, pages 1--22, 2021.

\bibitem{doikov2021minimizing}
Nikita Doikov and Yurii Nesterov.
\newblock Minimizing uniformly convex functions by cubic regularization of
  {N}ewton method.
\newblock {\em Journal of Optimization Theory and Applications}, pages 1--23,
  2021.

\bibitem{doikov2022affine}
Nikita Doikov and Yurii Nesterov.
\newblock Affine-invariant contracting-point methods for convex optimization.
\newblock {\em Mathematical Programming}, pages 1--23, 2022.

\bibitem{doikov2018randomized}
Nikita Doikov and Peter Richt{\'a}rik.
\newblock Randomized block cubic {N}ewton method.
\newblock In {\em International Conference on Machine Learning}, pages
  1289--1297, 2018.

\bibitem{dvurechensky2019near}
Pavel Dvurechensky, Alexander Gasnikov, Petr Ostroukhov, C{\'e}sar~A Uribe, and
  Anastasiya Ivanova.
\newblock Near-optimal tensor methods for minimizing the gradient norm of
  convex function.
\newblock {\em arXiv preprint arXiv:1912.03381}, 2019.

\bibitem{gasnikov2019near}
Alexander Gasnikov, Pavel Dvurechensky, Eduard Gorbunov, Evgeniya Vorontsova,
  Daniil Selikhanovych, C{\'e}sar~A. Uribe, Bo~Jiang, Haoyue Wang, Shuzhong
  Zhang, S{\'e}bastien Bubeck, Jiang Qijia, Yin~Tat Lee, Li~Yuanzhi, and
  Sidford Aaron.
\newblock Near optimal methods for minimizing convex functions with {L}ipschitz
  $p$-th derivatives.
\newblock In {\em Conference on Learning Theory}, pages 1392--1393, 2019.

\bibitem{goldfeld1966maximization}
Stephen~M. Goldfeld, Richard~E. Quandt, and Hale~F. Trotter.
\newblock Maximization by quadratic hill-climbing.
\newblock {\em Econometrica: Journal of the Econometric Society}, pages
  541--551, 1966.

\bibitem{grapiglia2017regularized}
Geovani~N. Grapiglia and Yurii Nesterov.
\newblock Regularized {N}ewton methods for minimizing functions with
  {H}\"{o}lder continuous {H}essians.
\newblock {\em SIAM Journal on Optimization}, 27(1):478--506, 2017.

\bibitem{grapiglia2019accelerated}
Geovani~N. Grapiglia and Yurii Nesterov.
\newblock Accelerated regularized {N}ewton methods for minimizing composite
  convex functions.
\newblock {\em SIAM Journal on Optimization}, 29(1):77--99, 2019.

\bibitem{grapiglia2020tensor}
Geovani~N. Grapiglia and Yurii Nesterov.
\newblock Tensor methods for minimizing convex functions with {H}\"{o}lder
  continuous higher-order derivatives.
\newblock {\em SIAM Journal on Optimization}, 30(4):2750--2779, 2020.

\bibitem{hanzely2020stochastic}
Filip Hanzely, Nikita Doikov, Peter Richt{\'a}rik, and Yurii Nesterov.
\newblock Stochastic subspace cubic {N}ewton method.
\newblock In {\em International Conference on Machine Learning}, pages
  4027--4038. PMLR, 2020.

\bibitem{kamzolov2020near}
Dmitry Kamzolov and Alexander Gasnikov.
\newblock Near-optimal hyperfast second-order method for convex optimization
  and its sliding.
\newblock {\em arXiv preprint arXiv:2002.09050}, 2020.

\bibitem{kantorovich1948functional}
Leonid~V. Kantorovich.
\newblock Functional analysis and applied mathematics. [in {R}ussian].
\newblock {\em Uspekhi Matematicheskikh Nauk}, 3(6):89--185, 1948.

\bibitem{kovalev2022first}
Dmitry Kovalev and Alexander Gasnikov.
\newblock The first optimal acceleration of high-order methods in smooth convex
  optimization.
\newblock {\em arXiv preprint arXiv:2205.09647}, 2022.

\bibitem{kovalev2019stochastic}
Dmitry Kovalev, Konstantin Mishchenko, and Peter Richt{\'a}rik.
\newblock Stochastic {N}ewton and cubic {N}ewton methods with simple local
  linear-quadratic rates.
\newblock {\em arXiv preprint arXiv:1912.01597}, 2019.

\bibitem{levenberg1944method}
Kenneth Levenberg.
\newblock A method for the solution of certain non-linear problems in least
  squares.
\newblock {\em Quarterly of applied mathematics}, 2(2):164--168, 1944.

\bibitem{lu2018relatively}
Haihao Lu, Robert~M. Freund, and Yurii Nesterov.
\newblock Relatively smooth convex optimization by first-order methods, and
  applications.
\newblock {\em SIAM Journal on Optimization}, 28(1):333--354, 2018.

\bibitem{marquardt1963algorithm}
Donald~W. Marquardt.
\newblock An algorithm for least-squares estimation of nonlinear parameters.
\newblock {\em Journal of the society for Industrial and Applied Mathematics},
  11(2):431--441, 1963.

\bibitem{mishchenko2021regularized}
Konstantin Mishchenko.
\newblock Regularized {N}ewton method with global ${O}(1/k^{2})$ convergence.
\newblock {\em arXiv preprint arXiv:2112.02089}, 2021.

\bibitem{monteiro2013accelerated}
Renato D.~C. Monteiro and Benar~F. Svaiter.
\newblock An accelerated hybrid proximal extragradient method for convex
  optimization and its implications to second-order methods.
\newblock {\em SIAM Journal on Optimization}, 23(2):1092--1125, 2013.

\bibitem{nesterov2008accelerating}
Yurii Nesterov.
\newblock Accelerating the cubic regularization of {N}ewton's method on convex
  problems.
\newblock {\em Mathematical Programming}, 112(1):159--181, 2008.

\bibitem{nesterov2018lectures}
Yurii Nesterov.
\newblock {\em Lectures on convex optimization}, volume 137.
\newblock Springer, 2018.

\bibitem{nesterov2019implementable}
Yurii Nesterov.
\newblock Implementable tensor methods in unconstrained convex optimization.
\newblock {\em Mathematical Programming}, pages 1--27, 2019.

\bibitem{nesterov2021inexact}
Yurii Nesterov.
\newblock Inexact accelerated high-order proximal-point methods.
\newblock {\em Mathematical Programming}, pages 1--26, 2021.

\bibitem{nesterov2021superfast}
Yurii Nesterov.
\newblock Superfast second-order methods for unconstrained convex optimization.
\newblock {\em Journal of Optimization Theory and Applications}, 191(1):1--30,
  2021.

\bibitem{nesterov2022quartic}
Yurii Nesterov.
\newblock Quartic regularity.
\newblock {\em arXiv preprint arXiv:2201.04852}, 2022.

\bibitem{nesterov1994interior}
Yurii Nesterov and Arkadi Nemirovski.
\newblock {\em Interior-point polynomial algorithms in convex programming}.
\newblock SIAM, 1994.

\bibitem{nesterov2006cubic}
Yurii Nesterov and Boris Polyak.
\newblock Cubic regularization of {N}ewton's method and its global performance.
\newblock {\em Mathematical Programming}, 108(1):177--205, 2006.

\bibitem{ortega2000iterative}
James~M. Ortega and Werner~C. Rheinboldt.
\newblock {\em Iterative solution of nonlinear equations in several variables}.
\newblock SIAM, 2000.

\bibitem{polyak2007newton}
Boris Polyak.
\newblock Newton’s method and its use in optimization.
\newblock {\em European Journal of Operational Research}, 181(3):1086--1096,
  2007.

\bibitem{polyak2009regularized}
Roman Polyak.
\newblock Regularized {N}ewton method for unconstrained convex optimization.
\newblock {\em Mathematical programming}, 120(1):125--145, 2009.

\bibitem{rodomanov2016superlinearly}
Anton Rodomanov and Dmitry Kropotov.
\newblock A superlinearly-convergent proximal {N}ewton-type method for the
  optimization of finite sums.
\newblock In {\em International Conference on Machine Learning}, pages
  2597--2605, 2016.

\bibitem{ueda2009regularized}
Kenji Ueda and Nobuo Yamashita.
\newblock A regularized {N}ewton method without line search for unconstrained
  optimization.
\newblock {\em Technical Report}, 2009.

\end{thebibliography}

\end{document}